\documentclass{amsart}		
\usepackage{verbatim,amsmath,latexsym,amssymb,amscd}

\theoremstyle{plain}
\newtheorem{thm}{Theorem}[section]

\newtheorem{lemma}[thm]{Lemma}
\newtheorem{prop}[thm]{Proposition}

\newtheorem{cor}[thm]{Corollary}

\theoremstyle{definition}
\newtheorem{definition}[thm]{Definition}
\newtheorem{example}[thm]{Example}

\def\mylabel#1{\label{#1}}

\input epsf

\begin{document}

\title[Polygon recognition and symmetry detection]{ Polygon recognition and symmetry detection}
\author{Mireille Boutin}
\today
\begin{abstract}
We introduce an approach based on moving frames for polygon recognition and symmetry detection. 
We  present detailed algorithms for recognition of polygons modulo the special Euclidean, Euclidean, equi-affine, skewed-affine and similarity Lie groups, and 
explain the procedure for a generic Lie group. 
The time complexity of our algorithms is linear in the number of vertices and they are noise resistant. The signatures  used allow the detection of partial as well as approximate equivalences. Our method is a particular case of a general method for curve recognition modulo Lie group action.
\end{abstract}
\maketitle

\section{Introduction}

This paper is devoted to equivalence of polygons under Lie group actions. As a subproblem, we also consider symmetries of polygons, which are nothing but self-equivalences,  i.~ e.~  non-trivial group transformations leaving the polygon unchanged.  
We are interested in  global, partial and approximate equivalences. The approach we suggest is based on the theory of moving frames (\cite{C}, \cite{FO1}, \cite{FO2})  and is a particular case of the method presented in \cite{MBjiscr} for  curve recognition.  This method consists in associating to every curve $C$ a certain polygon $P(C)$. 
The polygon is defined in such a way that if $\tilde{C}$ is another curve with $g\cdot\tilde{C}=C $ for some $g\in G$, then $g\cdot P(\tilde{C})=P(C)$. 
So as a first equivalence test, one can check  whether $P(\tilde{C})$ and $P(C)$ are equivalent under $G$. 
If not, then $\tilde{C}$ and $C$ are not equivalent under $G$. If  $P(\tilde{C})$ and $P(C)$ are equivalent under $G$, then we describe a detailed method to check whether the two curves are equivalent.

The main idea of our method for polygon recognition consists in constructing a joint invariant signature curve ($JIS$ curve) for every  polygon.
The signature curve of all polygons which belong to the same equivalence class is the same, and  symmetries of a polygon  manifest as  repetitions in its signature.

By construction, it is a very simple and visual approach, but more importantly it is general, in the sense that it has far more applicability than the particular cases we present here.
 In fact, it can be used for detecting equivalences under any Lie group which acts (locally) effectively on subsets, provided some slight regularity conditions.
Moreover, it can be generalized to higher dimensional structures such as polyhedra.  So although many algorithms for specific symmetry detection or equivalence in polygons are known ( see for example \cite{WWV}), we believe the algorithms we present here are interesting on their own, as they lay the basis of a general approach for equivalence of polygons and symmetry detection.

One advantage of this method is that it is noise resistant and can therefore be used for detection of approximate symmetries.  
Another advantage is that  each point  of the signature only depends on a few consecutive points of the  polygon. We are in fact able  to build  signatures which  indicate  partial equivalences, i.~ e.~ 
when two pieces of  a polygon are equivalent. Moreover, the dimension of the signature is optimal and  so is the complexity of the corresponding detection algorithms.

Our point is that the construction of a $JIS$ curve  is an easy, visual approach, which can be generalized to most Lie group actions on manifolds, and  the modern method of moving frame provides us with effective tools to compute the invariants we need.

In the following, we will  provide a full solution to the problem of detection of 
 all area preserving affine symmetries (rotations, reflections, equi-affine and skewed affine transformations). We will also provide a full solution for the problem of polygon recognition modulo the special Euclidean, full Euclidean, equi-affine, skewed-affine and similarity Lie groups. We will also explain our method  in details for a generic Lie group acting on a generic manifold.

\section{Mathematical Foundations}

Let $G$ be a Lie group acting  on a $m$-dimensional manifold $M$. 

\begin{definition} An {\em invariant} is a real valued function $I:M\rightarrow {\mathbb R}$ which remains unchanged under the action of $G$, more precisely
\[I(g\cdot p)=I(p), \text{ for all }p\in M \text{ and }g\in G.
\]
\end{definition}

\begin{definition}
We say that $G$ acts {\em freely} on $M$ if the identity is the only element of $G$ that fixes any point of $M$.
\end{definition}

\begin{definition}
 We say that $G$ acts {\em regularly} on $M$ if all orbits have the same dimension and if any point $p_0 \in M$  is surrounded by an arbitrarily small neighborhood whose intersection with the orbit through $p_0$ is connected.
\end{definition}

Most of our results are based on the following important theorem.
See \cite{Obook2} for a proof.

\begin{thm}[Frobenius Theorem]
\mylabel{fundamental theorem} If $G$ acts on an open set $O\subset M$ regularly with $s$ dimensional orbits, then  $\forall p_0 \in O$ there exist $m-s$ functionally independent invariants $I_1,\ldots,I_{m-s}$  defined on a neighborhood $U$ of $p_0$ such that any other invariant $H$ defined near $p_0$ is a function $H=f(I_1,\ldots,I_{m-s})$. Moreover, two points $p_1, p_2 \in U$ are in the same orbit if and only if $I_i (p_1)= I_i (p_2)$, for all $ i=1,\ldots ,m-s$.  \end{thm}

The set $\{ I_1,\ldots,I_{m-s} \}$ is often called a {\em complete fundamental set of invariants on $U$}.

 The modern theory of moving frames, as developed by Olver and Fels in \cite{FO1} and \cite{FO2} defines a (left) moving frame as follows.   

\begin{definition} A {\em moving frame} is a map $\rho : M\rightarrow G$ such that $\rho (g\cdot p)= g\cdot\rho(p)$,  $\forall p\in M$, $\forall g\in G$.

 A {\em local moving frame} is a map $\rho : M\rightarrow G$ such that $\rho (g\cdot p)= g\cdot\rho(p)$,  $\forall p\in M$, $\forall g\in N_e\subset G$, for some  neighborhood  $N_e$ of the identity $e\in G$.
\end{definition}

The conditions of existence of a moving frame are very precise.

\begin{thm} A (local) moving frame exists in a neighborhood of a point $p_0 \in M$ if and only if $G$ acts (locally) freely and regularly near $p_0$.
\end{thm}

\begin{definition}
We say that $G$ acts on $M$ effectively if
\[\{g\in G |\quad g\cdot p=p, \text{ for all }p\in M \}= \{e \}
\]

We say that $G$ acts on $M$ locally effectively if
\[\{g\in G |\quad g\cdot p=p, \text{ for all }p\in M \}
\]
is a discrete subset of $G$.
\end{definition}

\begin{definition}
We say that $G$ acts  effectively on subsets of $M$ if, for any open subset $U\subset M$,
\[\{g\in G |\quad g\cdot p=p, \text{ for all }p\in U \}= \{e \}.
\]

We say that $G$ acts locally effectively on subsets of $M$ if,  for any open subset $U\subset M$,
\[\{g\in G |\quad g\cdot p=p, \text{ for all }p\in M \}
\]
is a discrete subset of $G$.
\end{definition}

For analytic group actions, effectiveness implies effectiveness on subsets. However this is not true for general group actions.

A moving frame can be used as a tool to compute a complete fundamental set of invariants. See \cite{FO2} for a detailed algorithm. Let  $M^{\times (n)}:=M\times M \times \ldots \times M$ ($n$ times)  be the Cartesian product of $n$ copies of the manifold $M$. In the case where the action is not (locally) free, one option is to prolong the action of  $G$  on $M^{\times (n)}$    by  setting $g\cdot (p_1,\ldots ,p_k)=(g\cdot p_1,\ldots ,g\cdot p_k)$, for all $g\in G$ and  $(p_1,\ldots ,p_k)\in M^{\times (k)}$, and  hope that the action then becomes free.

Let $r$ be the dimension of $G$. The following  important result is proved in \cite{MBorbits}.

\begin{thm}
\mylabel{local freeness}
If $G$ acts (locally) effectively on subsets of $M$, then there exists a minimal integer $n_0$ such that, for all integers $n\geq n_0$, $G$ acts locally freely on an open and dense subset of  $M^{\times (n)}$.
\end{thm}

\begin{definition}
Let $n\in {\mathbb N}$. An {\em $n$-point joint invariant}, or  {\em joint invariant} for simplicity, is an invariant of the prolonged action of $G$ on $M^{\times (n)}$.
\end{definition}

Let  $H$ be a finite group acting on $k$ elements. In particular, $H$ acts on any given $k$ points  $p_1,\ldots ,p_k\in M$. So we have  an induced action of  $H$ on  $M^{\times (k)}$. Let $z^{(k)}\in M^{\times (k)}$. The following theorem will be used in this paper. Its proof is inspired by \cite{sturmfels} chapter two.

\begin{thm}
\mylabel{finite groups}
If $G$ acts regularly and  $H$ acts freely  in a neighborhood of $z^{(k)}$,  then in a   neighborhood $U$ of $z^{(k)}$ there exists a complete fundamental set of $G$-invariants $J_1,\ldots ,J_N: U\rightarrow {\mathbb R}$ which are also invariant under $H$.

Moreover,  we can choose $U$ such that two points $z^{(k)}_1\in h_1\cdot U$ and $z^{(k)}_2\in h_2\cdot U$ for $h_1,h_2\in H$ are in the same orbit relative to $G\times H$ if and only if $J_i(z^{(k)}_1)=J_i(z^{(k)}_2)$, for all $i=1,\ldots ,N$.
\end{thm}

\begin{proof}
By Theorem \ref{fundamental theorem}, there exist $\{ I_1,\ldots ,I_N \}$  a complete fundamental set of invariants under $G$ defined near  $z^{(k)}$. Define
\[P_i(t)=\Pi_{h \in H} (I_i(h \cdot (p_1,\ldots ,p_k))-t) ,\]
 for $i=1,\ldots ,N.$ We can view $P_i(t)$ as a polynomial in $t$ whose coefficients are functions of $p_1,\ldots ,p_k$. In fact, these coefficients are invariant under $G$. Moreover, since all $P_i(t)$'s are invariant under $H$,  their coefficients are also invariant under $H$.

Observe that $P_i(I_i(p_1,\ldots ,p_k) )=0$. In other words, there exists a non-trivial functional relationship between $I_i(p_1,\ldots ,p_k)$ and the coefficients of  $P_i(t)$. This means that, locally, $I_i(p_1,\ldots ,p_k)$ can be written as a function of the coefficients of  $P_i(t)$. Since $I_1,\ldots ,I_N$ are functionally independent, there must be $N$ functionally independent functions among the coefficients of the $P_i(t)$'s. This shows the first part of the statement.

To prove the second part, write  $z^{(k)}_j= h_j\cdot \bar{z}^{(k)}_j$, with  $\bar{z}^{(k)}_j\in U$ for $j=1,2$. By freeness of the action of $H$ in a neighborhood of $z^{(k)}$, we can choose $U$ so that $G\times H$ acts on $U$ regularly.
We have  $g\cdot z^{(k)}_1= z^{(k)}_2 $, for some $g\in G$ if and only if
$g\cdot h_1\cdot \bar{z}^{(k)}_1 =h_2\cdot \bar{z}^{(k)}_2$, or equivalently
$h_2^{-1}\cdot g\cdot h_1\cdot \bar{z}^{(k)}_1 = \bar{z}^{(k)}_2$.
By Theorem  \ref{fundamental theorem}, this happens if and only if  $J_i(\bar{z}^{(k)}_1)=J_i(\bar{z}^{(k)}_2)$, for all $i=1,\ldots ,N$. which  is equivalent to saying that
 $J_i( h_1\cdot  \bar{z}^{(k)}_1)=J_i(h_2 \cdot \bar{z}^{(k)}_2)$, for all  $i=1,\ldots ,N$ and the conclusion follows.
\end{proof}

\section{A signature for global polygon recognition}

\subsection{Equivalence of ordered sets of points under Lie group action}

In what follows, we will keep writing $M$ for a generic $m$ dimensional manifold and $G$ for a generic $r$ dimensional Lie group acting on $M$. Let  $(p_1,\ldots ,p_k)$ be a point of $M^{\times (k)}$. Suppose that $G$ acts on $M^{\times (k)}$ regularly with $s$-dimensional orbits. Then by theorem \ref{fundamental theorem}, there exist  fundamental invariants $I_1, \ldots ,I_{km-s}:U_M\subset M^{\times (k)}\rightarrow {\mathbb R} $. The map 
\begin{eqnarray}
S_M&:&U_M\subset M^{\times (k)}\rightarrow {\mathbb R}\nonumber \\
\text{defined by }\hspace{1cm}S_M(p_1,\ldots ,p_k)&=& \left(\begin{array}{c}
I_1(p_1,\ldots ,p_k)\\
\vdots\\
I_{mk-s}(p_1,\ldots ,p_k)
\end{array}
 \right) \nonumber
\end{eqnarray}
is a signature for ordered sets of $k$ points in the following sense.

\begin{thm}
\mylabel{equivalence ordered set of points}
Let $(p_1,\ldots ,p_k)$ and $(q_1,\ldots ,q_k)$ be two points of $U_M\subset M^{\times (k)}$.
There exists $g\in G$ such that $g \cdot (p_1,\ldots ,p_k)=(q_1,\ldots ,q_k)$ if and only if \[ S_M(p_1,\ldots,p_k)=S_M(q_1,\ldots,q_k).\]
\end{thm}

\begin{proof}
By theorem \ref{fundamental theorem}.
\end{proof}

So from the value of a finite set of invariants, one can conclude about the equivalence of two ordered sets of points. This provides an easy way to recognize ordered set of points up to Lie group action.

\subsection{Equivalence  of polygons under Lie group action}

Let the cyclic group ${\mathbb Z}_k$ act on  $M^{\times (k)}$ by  permuting the $k$ points cyclically. Let $\pi \in {\mathbb Z}_2$ act on $M^{\times (k)}$ by reversing the order of the $k$ points, i.~ e.~ $\pi (p_1,p_2,\ldots ,p_k)=(p_k,\ldots,p_2 ,p_1)$, for all $(p_1,p_2,\ldots ,p_k)\in M^{\times (k)}$.
Together, ${\mathbb Z}_k$ and ${\mathbb Z}_2$ generate a group acting on $M^{\times (k)}$. We shall call this group ${\mathbb H}_k= < {\mathbb Z}_k ,{\mathbb Z}_2  >$. 

\begin{lemma}
\mylabel{Hk group}
 If $h\in {\mathbb H}_k$, then either  
\begin{eqnarray}
 h&\in& {\mathbb Z}_k,\\
&\text{ or }&\nonumber \\
 h&=&c\cdot \pi , \text{ with } c\in {\mathbb Z}_k.
\end{eqnarray}

\end{lemma}


Let  $(p_1,\ldots ,p_k)$ be a point of $M^{\times (k)}$.
 Let ${\mathfrak P }^k=M^{\times (k)}/{\mathbb Z}_k$ be the set of $k$ ordered points in $M$  modulo the action of ${\mathbb Z}_k$ . If $p_1,\ldots ,p_k \in M$, the corresponding  point in ${\mathfrak P }^k $ will be written as $ \left[ p_1,\ldots ,p_k \right]  $. The action of $G$ on   $M^{\times (k)}$ naturally induces an action of $G$ on   ${\mathfrak P }^k$, namely $g\cdot \left[p_1,\ldots ,p_k \right]=\left[g\cdot p_1,\ldots ,g\cdot p_k \right]$, for  $g\in G$ and $p_1,\ldots ,p_k \in M$.

Let ${\mathcal P }^k=M^{\times (k)} / {\mathbb H}_k$ be the set of  $k$ ordered points in $M$  modulo the action of ${\mathbb H}_k$.   If $p_1,\ldots ,p_k \in M$, the corresponding  point in ${\mathcal P }^k $ will be written as $ \langle p_1,\ldots ,p_k \rangle  $. The action of $G$ on   $M^{\times (k)}$ naturally induces an action of $G$ on   ${\mathcal P }^k$, namely $g\cdot \langle p_1,\ldots ,p_k \rangle =\langle g\cdot p_1,\ldots ,g\cdot p_k \rangle$, for all  $g\in G$ and $p_1,\ldots ,p_k \in M$.

\begin{definition}
A {\em $k$-vertex polygon}, or {\em $k$-gon},  is a point of  ${\mathcal P }^k$.
\end{definition}

\begin{definition}
We say that two $k$-gons $P= \langle p_1,\ldots ,p_k \rangle $ and $Q= \langle q_1,\ldots ,q_k \rangle $ are {\em equivalent under $G$} if there exist $g\in G$ such that $g\cdot \langle p_1,\ldots ,p_k \rangle =\langle q_1,\ldots ,q_k \rangle $. In that case, we write $P\equiv  Q \mod G$. 
\end{definition}

\begin{definition}
We say that a polygon  $P= \langle p_1,\ldots ,p_k \rangle $ {\em has a $G$-symmetry} if there exists $g\in G \setminus \{ e\} $ such that $g\cdot\langle  p_1,\ldots ,p_k \rangle =\langle  p_1,\ldots ,p_k \rangle $.  
\end{definition}

Suppose that a polygon $P=\langle p_1,\ldots ,p_k \rangle $ has a $G$-symmetry. This means that there exists $g\in G\setminus \{e \}$ and $h\in{\mathbb H}_k $ such that $g\cdot( p_1,\ldots ,p_k )=h\cdot ( p_1,\ldots ,p_k )$. 
According to Lemma \ref{Hk group}, either  $h\in {\mathbb Z}_k$ or $h=c\cdot \pi$, with $c\in {\mathbb Z}_k$ and $\pi \in {\mathbb Z}_2$ as defined above.
Similarly if  $P=\{p_1,\ldots ,p_k \}$ is equivalent to $Q=\{q_1,\ldots ,q_k \}$ modulo $G$, then $g\cdot (p_1,\ldots ,p_k)=h\cdot (q_1,\ldots ,q_k)$ for some $g\in G$ and $h\in {\mathbb Z}_k$ or $h=c \pi$ with $c\in {\mathbb Z}_k $. These facts will be use to simplify our symmetry detection and polygon recognition algorithms later on.

If $G$ acts regularly with $s$-dimensional orbits on some open set $U\subset M^{\times (k)}$,  
then by Theorem \ref{fundamental theorem} , in a neighborhood $U_M$ of any point $(p_1,\ldots ,p_k)\in U  $,   there exists a complete fundamental set of $G$-invariants $\{I_1, \ldots , I_{mk-s}\}:U_1 \rightarrow {\mathbb R}$.
Assuming that $H_k$ acts freely on $U$ (which can be guaranteed  by taking $p_1,\ldots ,p_k$ distinct for example, and by choosing $U$ small enough),   
then by Theorem \ref{finite groups},  there also exists  $\{ \bar{I}_1,\ldots , \bar{I}_{mk-s} \}$, a complete fundamental set of $(G \times {\mathbb Z}_k ) $-invariants, as well as  $\{ \tilde{I}_1,\ldots , \tilde{I}_{mk-s} \}$, a complete fundamental set of $(G \times {\mathbb H}_k ) $-invariants, all defined on some neighborhood of $(p_1,\ldots ,p_k)$.     
 The maps

\begin{eqnarray}
S_{\mathfrak M} &:&U_{\mathfrak M}\subset{\mathfrak M}^k\rightarrow {\mathbb R}\nonumber \\
S_{\mathfrak M}(\left[ p_1,\ldots ,p_k \right] )&=& \left(\begin{array}{c}
\bar{I}_1(p_1,\ldots ,p_k)\\
\vdots\\
\bar{I}_{mk-s}(p_1,\ldots ,p_k)
\end{array}
 \right), \nonumber
\end{eqnarray}

\begin{eqnarray}
S_{\mathcal P} &:& U_{\mathcal P}\subset {\mathcal P}^k\rightarrow {\mathbb R}\nonumber \\
S_{\mathcal P}(\{ p_1,\ldots ,p_k\} )&=& \left(\begin{array}{c}
\tilde{I}_1(p_1,\ldots ,p_k)\\
\vdots\\
\tilde{I}_{mk-s}(p_1,\ldots ,p_k)
\end{array}
 \right) \nonumber
\end{eqnarray}

and $S_M$
 constitute signatures for polygons in the following sense.

\begin{thm}
\mylabel{3 signature polygons}
Let $P=\langle p_1 ,\ldots ,p_k\rangle $ and   $Q=\langle q_1 ,\ldots ,q_k\rangle $ be two $k$-gons. Assume that $(p_1,\ldots ,p_k) \in U_M$, $\left[ p_1,\ldots ,p_k \right] \in U_{\mathfrak M}$ and that $\{ p_1,\ldots ,p_k\} \in U_{\mathcal P} $.  Then $P\equiv Q \mod G$ 
\begin{eqnarray}
\Leftrightarrow    S_M( p_1,\ldots ,p_k ) & = & S_M(h\cdot (q_1,\ldots ,q_k)) \text{ for some } h\in {\mathbb H}_k\nonumber \\
\Leftrightarrow   S_{\mathfrak M}( \left[p_1,\ldots ,p_k \right]) & = & S_{\mathfrak M}( \left[ q_1,\ldots ,q_k \right] ) \nonumber\\
& \text{ or }& \nonumber \\  
S_{\mathfrak M}( \left[p_1,p_2, \ldots ,p_k \right]) & = & S_{\mathfrak M}( \left[ q_k,\ldots, q_2 ,q_1 \right] ) \nonumber\\
\Leftrightarrow   S_{\mathcal P}( \langle p_1,\ldots ,p_k \rangle )& = & S_{\mathcal P}(\langle q_1,\ldots ,q_k \rangle )\nonumber
\end{eqnarray}
\end{thm}

\begin{proof}
The necessity of the first and third statements follow from the invariance of the signature. For the second statement, we also use lemma \ref{Hk group}.

To prove sufficiency of the first statement, assume that $S_1(p_1,\ldots ,p_k)=S_1(h\cdot (q_1,\ldots ,q_k))$. Then by Theorem \ref{finite groups}, there exists $g\in G$ such that $g\cdot (p_1,\ldots ,p_k)=h\cdot (q_1,\ldots ,q_k) $. So $g\cdot \{p_1 , \ldots ,p_k \}= \{  q_1,\ldots ,q_k \}$.

To prove sufficiency of the second statement, assume that
 $S_2 (\left[ p_1,\ldots ,p_k \right] )=S_2 ( \left[ z\cdot ( q_1,\ldots ,  q_k) \right] )$
 for some $z\in {\mathbb Z}_2$.
 By theorem \ref{finite groups}, this means that there exists $g\in G$ and $c\in {\mathbb Z}_k$ such that 
\[
g\cdot (p_1,\ldots ,p_k) =c\cdot (z\cdot q_1, \ldots ,z\cdot q_k),
\]
and therefore $g\cdot \{p_1,\ldots ,p_k  \}=\{q_1,\ldots ,q_k \}$.
The proof for sufficiency of the third statement is similar.
\end{proof}

\subsection{Equivalence of point configurations under Lie group action}

Let ${\mathbb S}_k$ be the symmetric group. The elements of  ${\mathbb S}_k$ act on $\{ (p_1,\ldots ,p_k )\in M^{\times (k)}\}$ by permuting the points $p_1,\ldots ,p_k$. 

\begin{definition} A $k$-point configuration is a point of  $M^{\times (k)}/{\mathbb S}_k$.
\end{definition}

In other words, a $k$-point configuration on $M$ is  a finite set of $k$ points on $M$ which are not ordered in any way. 
We shall use the notation $ | p_1,\ldots ,p_k| $ for the  $k$-point configurations corresponding to $p_1,\ldots ,p_k \in M$.
  The action of $G$ on $M$ naturally induces an action of $G$ on  $k$-point configurations. 
The polygon recognition method is easy to extend for point configuration recognition. 
In fact, we can repeat the same arguments as before to claim the existence of a complete set of fundamental invariants $\hat{I_1},\ldots ,\hat{I}_{mk-s}$
under $G \times {\mathbb S}_k $.
 Also, the map

\begin{eqnarray}
S_4&:&M^{\times (k)}/{\mathbb S}_k \rightarrow {\mathbb R}\nonumber \\
S_4( | p_1,\ldots ,p_k| )&=& \left(\begin{array}{c}
\hat{I}_1(p_1,\ldots ,p_k)\\
\vdots\\
\hat{I}_{mk-s}(p_1,\ldots ,p_k)
\end{array}
 \right) \nonumber
\end{eqnarray}

is a  signature for $k$-point configurations  in the following sense.

\begin{thm}
Let $|p_1 ,\ldots ,p_k|$ and   $|q_1 ,\ldots ,q_k|$ be two $k$-point configurations. Then there exists $g\in G$ such that $g\cdot |p_1 ,\ldots ,p_k|= |q_1 ,\ldots ,q_k|$ 
\[ \Leftrightarrow S_4(| p_1,\ldots ,p_k |)= S_4(|q_1,\ldots ,q_k|).\]
\end{thm}

We can even go further. For example, we can consider a finite number of polygons (without any order). In a similar manner, we can define a signature which will characterize  these polygons up to Lie group action. This can be used in the case we want to recognize pictures made  of a finite number of disconnect pieces, each  being a polygon.

\section{A signature for partial polygon recognition}

The previous sections provide us with a way to recognize polygons globally. However, we are also interested in the case where a piece of a polygon is equivalent to a piece of another polygon. In particular, we would like the signature to indicate whether  pieces of two polygons are equivalent under group action, or if a polygon has a certain symmetry. This would be a complex task using the previous signatures. In the following, we  explain a simpler method.

Recall that $m$ is the dimension of the manifold $M$. As mentioned before, we would like to parameterize a signature with no more than $m$ invariants, since this is the optimal number. We will explain shortly how to choose suitable invariants. But first let us give some definitions.

 Given $m$ invariants $I_1,\ldots ,I_m:  M^{\times (n)}\rightarrow {\mathbb R}$ of the action of $G$ on $M^{\times (n)}$ and $(p_1,\ldots ,p_k)\in M^{\times (k)}$, define 
\[ I_{i,r}: M^{\times (k)}\rightarrow {\mathbb R}, \text{ for }i=1,\ldots m \text{ and } r=1,\ldots ,k \]\[\text{by }I_{i,r}(p_1, \ldots ,p_k  )= I_i (p_r, p_{r+1}, \ldots ,p_{r+n-1}). \] 
setting $p_{k+1}=p_1,\quad p_{k+2}=p_2,\quad \ldots \quad ,p_{k+n-1}=p_{n-1}$.

Define $S: M^{\times (k)}\rightarrow ({\mathbb R}^m)^{\times (k)}$ by 
\[ S(p_1,\ldots ,p_k)= \left( \begin{array}{c c c}
I_{1,1}(p_1,\ldots ,p_k)&,\ldots ,&I_{1,k}(p_1,\ldots ,p_k)\\
\vdots& & \vdots \\
I_{m,1}(p_1,\ldots ,p_k)&,\ldots ,&I_{m,k}(p_1,\ldots ,p_k)   
\end{array}
\right)
 \]

If we let  ${\mathfrak S}$ be the map  ${\mathfrak S}: {\mathfrak P }^k\rightarrow ({\mathbb R}^m)^{\times (k)} \text{ mod }{\mathbb Z}_k$   given by

\[{\mathfrak S}(\left[p_1,\ldots ,p_k   \right] )= S(p_1,\ldots ,p_k) \text{ mod }{\mathbb Z}_k,  
\]
 then the following diagram commutes.
\[
\begin{CD}
M^{\times (k)}     @>{S}>> ({\mathbb R}^m)^{\times (k)}\\
@V{\mod {\mathbb Z}_k}VV          @VV{\mod {\mathbb Z}_k }V\\
{\mathfrak P }^k         @>>{{\mathfrak S}  }> ({\mathbb R}^m)^{\times (k)}\mod {\mathbb Z}_k 
\end{CD}
\]

The maps $S$ and ${\mathfrak S} $ can be used as signatures in the following instances.

\begin{thm}[For global recognition]
\mylabel{recognition with m invariants}
Let $P= \langle p_1, \ldots ,p_k   \rangle $ and $Q=\langle q_1, \ldots ,q_k   \rangle$ be two $k$-gons. Assume that the  set    $\{I_{1,r}, \ldots ,I_{m,r}  \}_{r=1}^k$ contains a complete fundamental set of $k$-point joint invariants on some open set $U_k\subset M^{\times (k)}$ and that $ (p_1, \ldots ,p_k),(q_1, \ldots ,q_k)\in U_k$. Then  $P\equiv  Q$ modulo $G$
\begin{eqnarray}
\Leftrightarrow &  S(p_1,\ldots ,p_k)&=S(h \cdot (q_1,\ldots,q_k)), \text{ for some } h\in   {\mathbb H}_k,\nonumber\\
\Leftrightarrow&  {\mathfrak S}(\left[ p_1, \ldots ,p_k   \right] )&={\mathfrak S}( \left[ q_1, \ldots , q_k   \right] )  \nonumber \\
&&\text{or} \nonumber \\ 
&{\mathfrak S}(\left[ p_1, \ldots ,p_k   \right] )&={\mathfrak S}( \left[ q_k, \ldots , q_1)   \right] )\nonumber
\end{eqnarray}

\end{thm}

\begin{proof}
Follows from theorem \ref{3 signature polygons}.
\end{proof}

\begin{cor}[For symmetry detection.]
\mylabel{symmetry detection with m invariants}
Let $P= \langle p_1, \ldots ,p_k   \rangle $ be a polygon. Assume that  the  invariants    $\{I_{1,r}, \ldots ,I_{m,r}  \}_{r=1}^k$ contain a complete fundamental set of $k$-point joint invariants on some open set $U_k\subset M^{\times (k)}$ and that ${\mathbb H}_k\cdot (p_1, \ldots ,p_k)\subset U_k$. Then  $P$ has a  $G$-symmetry if and only if
\[ 
\Leftrightarrow  S(p_1,\ldots ,p_k)=S(h \cdot (p_1,\ldots,p_k)), \text{ for some } h\in   {\mathbb H}_k\setminus \{ e \}.
\]
\end{cor}

We now need to explain how to construct $m$ suitable invariants $I_1,\ldots ,I_m$, suitable in the sense that  $\{I_{1,r},\ldots , I_{m,r}\}_{r=1}^k $ contains  a complete fundamental set of $k$-point joint invariants on some open set.
Before we present the general method for polygon recognition, let us consider two instructive examples.

\section{An example of orientation preserving Lie group action}

\subsection{Construction of the signature}

As a first example, consider $SE(2)$, the group of orientation preserving rigid motion in the plane ( i.~ e.~   rotations and translations.) We call it the {\em special Euclidean group}. It is generally accepted to call the corresponding symmetries of polygons {\em rotational symmetries}, since  any such symmetry corresponds to  a rotation around some interior point of the polygon. 

This  well known result can be proved  using the moving frame method (see \cite{Ojis}.)

\begin{thm}For  $SE(2)^\curvearrowright {\mathbb R}^2$, we have the following. 

\begin{enumerate}

\item There are no one-point joint invariants.

\item There is one fundamental two-point joint invariants $I(p_1,p_2):{\mathbb R}^2\times {\mathbb R}^2\rightarrow {\mathbb R}$, namely $I(p_1,p_2)=|p_2-p_1|$.

\item There are three fundamental three-point joint invariants 
\[
I_1(p_1,p_2,p_3),\quad I_2(p_1,p_2,p_3),\quad I_3(p_1,p_2,p_3):{\mathbb R}^2\times {\mathbb R}^2 \times {\mathbb R}^2 \rightarrow {\mathbb R} ,
\]
namely
\begin{eqnarray}
 I_1(p_1,p_2,p_3)&=&|p_2-p_1|,\nonumber \\
  I_2(p_1,p_2,p_3)&=&|p_3-p_2|,\nonumber \\
 I_3(p_1,p_2,p_3)&= &\frac{1}{2}\det[z_{3}-z_{1},z_{2}-z_{1}] ,\text{ the signed area of the triangle}.\nonumber \\
&=: &\Delta_{123}.\nonumber
\end{eqnarray}

\end{enumerate}
\end{thm}

We are looking for two suitable joint invariants $J_1$ and $J_2$ to build a signature. Again by suitable we mean that $\{J_{1,r},J_{2,r} \}_{r=1}^k$ should contain a complete set of fundamental $k$-point joint invariants on some open set.

\begin{figure}[here]
\caption{}
\label{triangle}
\vspace{0.5cm}
\centerline{
\hbox{
\epsfysize=4.0cm
\epsfxsize=5.0cm
{\leavevmode \epsfbox{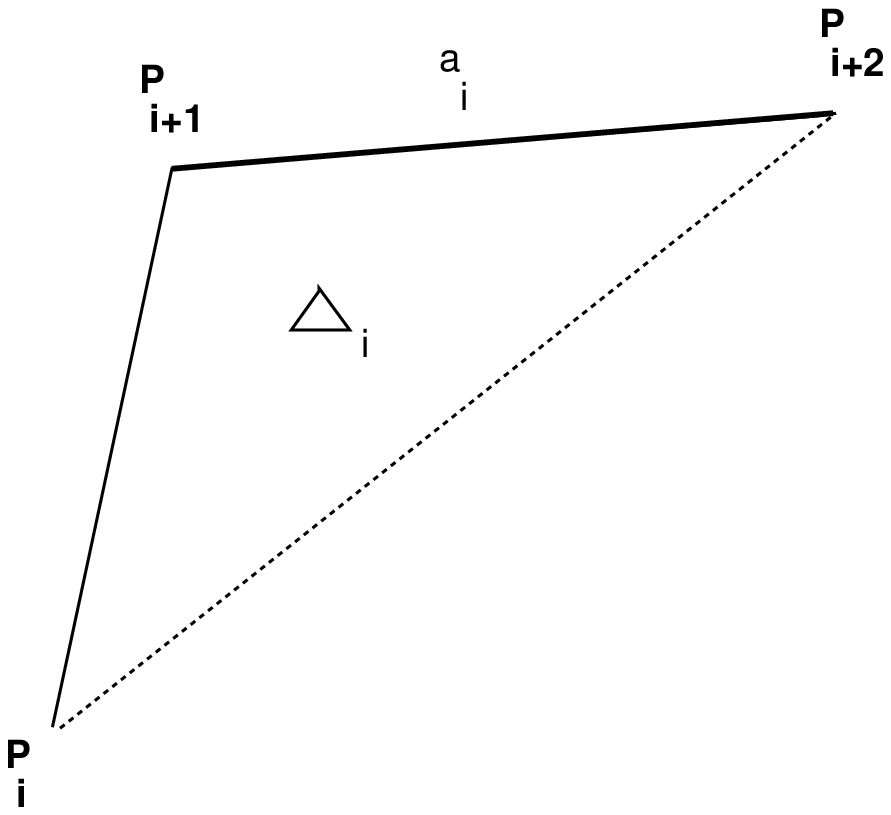}}
}
}
\end{figure}

In this case, we can take $J_1(p_1,p_2,p_3)=|p_3-p_2|$ and  $J_2(p_1,p_2,p_3)=\Delta_{123}$. According to our general method to be explained later (see Theorem \ref{star and starstar exist}), this is a natural choice.
Given a $k$-gon $P=\langle p_1,\ldots ,p_k \rangle $, we define its  special Euclidean joint invariant signature ($SEJIS$)  as the sequence of $k$ points given by

\begin{equation}
SEJIS(p_1,\ldots ,p_k)=\{a_i, \Delta_i   \}_{i=1}^k,
\end{equation} 

where  $a_i=|p_{i+2}-p_{i+1}|$ and $\Delta_i$ is the signed area given by the  determinant of the $2\times 2$ matrix $\frac{1}{2}[p_{i+2}-p_i,p_{i+1}-p_{i}]$. See the illustration in  figure \ref{triangle}.
The reasons why the invariants $a$ and $\Delta$ can be used to build  a signature are contained in two properties.

First, when evaluating $a_i$ and $\Delta_i$ for $i=k-1,k$ and $1$,  one obtains all fundamental joint invariants which only depend on the first two points $p_1$ and $p_2$. In this case, there is only one, namely $I(p_1,p_2)=|p_2-p_1|$. In other words, we have
\[\left\{ |p_2-p_1| \right\}\subset \{a_i,\Delta_i \}_{i=k-1}^1 \]
with $\{ |p_2-p_1| \}$  a complete fundamental set of  joint invariants only depending on $p_1$ and $p_2$. 
This guarantees the first property  called {\em two-point projectability}  ($\star$).  
Observe that if $I(p_1,p_2)=I(q_1,q_2)$, then there exists $ g\in SE(2)$ such that $g\cdot q_1=q_1$ and $g\cdot p_2=q_2$. This is the key idea in ($\star$).

Secondly, given $p_1$ and $p_2$ with $p_1\not= p_2$, then $p_3$ is uniquely determined by the value of $a_1=|p_3-p_2|$ and $\Delta_1=\frac{1}{2}\det [p_3-p_1,p_2-p_1]$. In other words, $p_3$ is a function 
\[ p_3=f(p_1,p_2,a_1,\Delta_1),\]
provided that $p_1\not= p_2$. In fact,   $p_{i+2}$ is a function $p_{i+2}=f(p_i,p_{i+1},a_i,\Delta_i)$ for all $i$'s, whenever $p_{i+1}\not= p_i$. This guarantees the second  property  called {\em third point reductivity} (property ($\star\star$)) when consecutive points are distinct.

As will be proved in Theorem \ref{starstar does not exist}, $n=3$ is the minimal number of points for which we can find two $n$-point joint invariants  $I_1$ and $I_2$ which are $(n-1)$-point projectable and $n^{th}$ point reductive on some open set.
 There are of course other suitable choices of invariants than $a_i$ and $\Delta_i$. As we will see,  ($\star$) and ($\star\star$) are enough to guarantee that $S=\{ J_{1,r}, J_{2,r} \}_{r=1}^k$ contains a complete fundamental set of $k$ point joint invariants and can therefore be used as a signature.

\begin{definition}
Choose an orientation (i.~ e.~ a traveling direction on the vertices) for $P=\langle p_1,\ldots ,p_k \rangle$.
The  {\em  special Euclidean joint invariant signature curve} ($SEJIS$ curve) of $P$ with respect to this orientation is the piecewise linear curve obtained by joining the points of the signature which correspond to consecutive vertices of the polygon by a straight oriented line.
\end{definition}

The $SEJIS$ curve represents the special Euclidean signature up to cyclic permutations of its $k$ points. This takes care of the ambiguity about the starting point $p_1$.
There remains one ambiguity: the traveling direction. In fact, the points of the $SEJIS$  are not invariants under the action of reversing the order of the vertices of the polygon.
 However if we restrict ourselves to simple polygons, i.~ e.~ polygons whose edges do not cross each other, then we can prescribe a specific orientation (clockwise for example) and this orientation remains unchanged under the action of $SE(2)$. In fact the $SEJIS$ curve characterizes all simple polygons.

\begin{thm}[For simple polygon recognition modulo $SE(2)$]
Two planar polygons whose edges do not cross and whose points are labeled clockwise are equivalent under the action of $SE(2)$ if and only if their $SEJIS$ curve with respect to the clockwise orientation is the same.
\end{thm}

\begin{proof}
Since the points of the signatures are functions of the basic $SE(2)$-invariants, they are $SE(2)$-invariant themselves. Moreover the order of the points is chosen in an invariant way, except for the starting point. Therefore if two polygons are equivalent under the action of $SE(2)$, then their signature will be identical, up to cyclic permutation. Now suppose that $P=\langle p_1,\ldots ,p_k  \rangle $ and $Q=\langle q_1,\ldots ,q_k  \rangle $ are two polygons with the same  $SEJIS=(s_1,\ldots , s_k   )$. 
Assume that $s_1$ corresponds to $(p_1, p_2, p_3)$ and  $(q_1, q_2, q_3)$, that    $s_2$ corresponds to $(p_2, p_3, p_4)$ and  $(q_2, q_3, q_4)$, and so on. 
Since the signature of the two polygons is the same,  we have $|p_2-p_1|=|q_2-q_1|$  (by  ($\star$)) . So we can find $g\in SE(2)$ which maps $p_1$ to $q_1$ and $p_2$ to $q_2$.
Moreover since   $p_{i+2}$ is uniquely prescribed by $p_i, p_{i+1}$ and the value of $a_i$ and $\Delta_i$ (by ($\star\star$)), we have that $g$ also maps $p_3$ to $q_3$, and $p_4$ to $q_4$, and so on. Therefore  $g\cdot P=Q$. 
\end{proof}

However, we do not need to restrict ourselves to simple polygons. All we have to do in order to characterize all polygons is to use our very same $SEJIS$ curve while taking into account the fact that we might have chosen  a different orientation and starting point.

\begin{thm}[For polygon recognition modulo $SE(2)$]
Two planar polygons $P=\{p_1,\ldots ,p_k \}$ and $Q=\{q_1,\ldots ,q_k \}$ are equivalent under the action of $SE(2)$ 
\[
\Leftrightarrow   SEJIS(p_1,\ldots ,p_k) =SEJIS(h\cdot (q_1,\ldots ,q_k)), \text{ for some } h\in{\mathbb H}_k.\]

\end{thm}

Unfortunately, the fact that we only characterize polygons up to ${\mathbb H}_k$ is inherent to the construction of the signature. However,  Lemma \ref{Hk group} facilitates the search for a possible $h \in {\mathbb H}_k$. In facts, since the $SEJIS$ commutes with rotations, we have the following useful lemma.

\begin{lemma}
\mylabel{simplify search for h in SE(2)}
Let 
\begin{eqnarray*}
SEJIS(p_1,\ldots ,p_k)&=&(s_1,\ldots ,s_k),\\
SEJIS(q_1,\ldots ,q_k)&=&(\sigma_1,\ldots ,\sigma_k),\\
\text{and }\quad   SEJIS(q_k,\ldots ,q_1)&=&(\bar{\sigma}_k,\ldots ,\bar{\sigma}_1).
\end{eqnarray*}
Then  $SEJIS(p_1,\ldots ,p_k)=SEJIS(h\cdot (q_1,\ldots ,q_k))$ for some $h\in {\mathbb H}_k$ if and only if
\begin{eqnarray}
 (s_1,\ldots ,s_k)&=&c\cdot (\sigma_1,\ldots ,\sigma_k)\nonumber \\
&\text{ or } & \nonumber \\
(s_1,\ldots ,s_k)  &= &c\cdot (\bar{\sigma}_k,\ldots ,\bar{\sigma}_1) \nonumber
\end{eqnarray}
for some $c\in {\mathbb Z}_k$.
\end{lemma}

Since a symmetry is a self-equivalence, we can modify the previous theorem in order to detect symmetries. In fact in this case, the ambiguity about the direction is waived and  the orientation of the polygon can be chosen arbitrarily.

\begin{thm}[For $SE(2)$-symmetry detection]
If $P=\langle p_1,\ldots,p_k \rangle$ is a planar polygon and $SEJIS(p_1,\ldots ,p_k)=(s_1,\ldots ,s_k)$, then $P$ has an $f$-fold rotational symmetry if and only if  
\[(s_{\frac{k}{f}+1}, \ldots, s_k,s_1,\ldots ,s_{\frac{k}{f}}) = (s_1,\ldots ,s_k),\]
 in other words, if and only if 
   the signature curve winds $f$ times on itself.
\end{thm}

\begin{proof}

The polygon  $P$ has an $f$-fold symmetry if and only if there exists $g\in SE(2)$ such that
$g\cdot (p_1,\ldots ,p_k)=(p_{\frac{k}{f}+1}, \ldots ,p_k, p_1, \ldots ,p_{\frac{k}{f}})$.
By invariance of the signature, that means 
\begin{eqnarray}
SEJIS(p_1,\ldots ,p_k)&=&SEJIS(p_{\frac{k}{f}+1}, \ldots ,p_k, p_1, \ldots ,p_{\frac{k}{f}})\nonumber \\
\Leftrightarrow (s_1,\ldots ,s_k)&=&(s_{\frac{k}{f}+1},\ldots ,s_{k}, s_1, \ldots , s_{\frac{k}{f}}),\nonumber
\end{eqnarray}
 which proves the necessity of the statement.

Now if $(s_1,\ldots ,s_k)=(s_{\frac{k}{f}+1},\ldots ,s_{k}, s_1, \ldots , s_{\frac{k}{f}})$, then by property ($\star$) and ($\star\star$), there exists $g\in SE(2)$ such that $g\cdot (p_1,\ldots ,p_k)= (p_{\frac{k}{f}+1},...,p_{k},p_{k+1},\ldots ,p_{\frac{k}{f}}) $. This proves sufficiency.
\end{proof}
 
So $P$ has an  $f$-fold symmetry if and only if the signature curve is traced $f$ times in the same direction as one travels along the curve. 
This can be checked in $O(k)$ by a computer. 



We implemented the algorithm using Matlab and computed the results for a few examples. One of them is   the four-branch star of figure \ref{star}. For a counterclockwise orientation, the program gives the following $SEJIS$ (rounded to 4  digits) for this eight-vertex polygon.

\[
SEJIS=\left[ 
\begin{array}{cccccccc} 
2.236 & 2.236 & 2.236 & 2.236 & 2.236 & 2.236 & 2.236 & 2.236  \\
4 & -3 & 4 & -3 & 4 & -3 & 4 & -3 
\end{array}
\right]
\]

The $SEJIS$ curve, represented in figure \ref{star},  is obtained by joining those points with a straight oriented segment. Although the polygon has eight vertices, the graph of the signature shows only two vertices: the signature curve winds four times on itself. This reflects the fact the four-branch star shown has a four-fold symmetry.

\subsection{Advantages of this $SEJIS$}

This is clearly  not the only way to build a signature. So why do we prefer this method to others?

First of all, this signature will indicate whether two pieces of polygons are the same (partial equivalences) up to Lie  group action. (See definition below.) This is because the invariants used depend on very few points and they are chosen so that partial equivalences correspond to  specific  similarities of the signature curves. More precisely, we have the following theorem.

\begin{thm}[For partial equivalences modulo $SE(2)$]
Let $P=\langle p_1,\ldots ,p_k\rangle$ be a planar polygon and $SEJIS(p_1,\ldots,p_k)=(s_1,\ldots ,s_k)$.
Let $Q=\langle q_1,\ldots ,q_l \rangle$ be another planar polygon  and $SEJIS(q_1,\ldots,q_k)=(\sigma_1,\ldots ,\sigma_l)$. Let $n\in {\mathbb N}$, let  $\tilde{P}=(p_{i},\ldots ,p_{i+n})$  and let $\tilde{Q}=(q_{j},\ldots ,q_{j+n})$. Let $s_i^x$ be the first component of $s_i$, namely $|p_{i+2}-p_{i+1}|$, and similarly for $\sigma_i^x$.
There exists $g\in SE(2)$ such that $g\cdot(p_{i} ,p_{i+1})=(q_{j} ,q_{j+1}) $ if and only if $s^x_{i-1}=s^x_{j-1}$.

For $n >1$,
there exist $g\in SE(2)$ such that $g\cdot(p_{i},\ldots ,p_{i+n})=(q_{j},\ldots ,q_{j+n}) $ if and only if 
$(s_{i}, \ldots ,s_{i+n-2})=(\sigma_{j}, \ldots ,\sigma_{j+n-2})$ and $s^x_{i-1}=s^x_{j-1}$.

\end{thm}

We call the $(n+1)$ consecutive vertices of $P$ given by $\tilde{P}=(p_{i},\ldots ,o_{i+n})$ a {\em piece of $P$}. If a piece of $P$ is equivalent to a piece of $Q$, more precisely  if 
\begin{eqnarray*}
(p_i,\ldots ,p_{i+n})&\equiv& (q_k,q_{k+1},\ldots ,q_{j+n})\mod G\\
\text{or }\quad (p_i,\ldots ,p_{i+n})&\equiv &(q_{j+n},\ldots , q_{k+1} ,q_{k})\mod G,
\end{eqnarray*}   
then we say that $P$ is {\em partially equivalent}  to $Q$.






From then, it is easy to modify our method in order to recognize what we call {\em polygonal segments} (or {\em open polygons}). Given an ordered set of points $(p_1,\ldots ,p_k)$ in the plane, define its signature as the ordered set of points given by

\begin{cor}[For open polygon recognition]
There exists $g\in SE(2)$ such that $g\cdot(p_1,\ldots ,p_k)=(q_1,\ldots ,q_k) $ if and only if 
\[ SEJIS_{ open}(p_1,\ldots ,p_k)=SEJIS_{open}(q_1,\ldots ,q_k). \] 
\end{cor}

Another advantage of our signature is that it is noise resistant. This is because it is parameterized by  functions which do not depend on derivatives. We are using joint invariants and in general, the value of such invariants does not change much when the points are slightly perturbed. In fact in this specific case, if we measure the points $p_i$ as $\tilde{p}_i$, and if the noise is such that the  measures  are within a certain error say
\[\tilde{p}_i=p_i \pm (\epsilon, \epsilon),
\] 
then the measured signature $\{\tilde{a}_i,\tilde{\Delta}_i \}_{i=1}^k$ will have the following precision:
\begin{eqnarray}
\tilde{a}_i&=&a_i\pm2\sqrt{2}\epsilon\nonumber \\
\tilde{\Delta}_i&=& \Delta_i \pm \frac{13}{2}(|p_{i+1}-p_{i}|+|p_{i+2}-p_{i}|+|p_{i+2}-p_{i+1}|)       \sqrt{\epsilon}\nonumber
\end{eqnarray}

\section{An example of non-orientation preserving Lie group action}

A slightly more complicated case is the recognition of planar polygons up to rotations and reflections. The corresponding group  is called the (full) Euclidean group $E(2)$ and consists in all rigid motion in the plane: translations, rotations and reflections. It is  not orientation preserving.

The following theorem can be obtained using the moving frame  method.

\begin{thm}
\mylabel{123 pt E(2)}
For  $E(2)^\curvearrowright {\mathbb R}^2$, the situation is as follows. 
\begin{enumerate}

\item There are no one-point joint invariants.

\item There is one fundamental two-point joint invariants $I(p_1,p_2):{\mathbb R}^2\times {\mathbb R}^2\rightarrow {\mathbb R}$, namely $I(p_1,p_2)=|p_2-p_1|$.

\item There are three fundamental three-point joint invariants 
\[
I_1(p_1,p_2,p_3),\quad I_2(p_1,p_2,p_3),\quad I_3(p_1,p_2,p_3):{\mathbb R}^2\times {\mathbb R}^2 \times {\mathbb R}^2 \rightarrow {\mathbb R} ,
\]
namely
\begin{eqnarray}
 I_1(p_1,p_2,p_3)&=&|p_2-p_1|,\nonumber \\
  I_2(p_1,p_2,p_3)&=&|p_3-p_2|,\nonumber \\
 I_3(p_1,p_2,p_3)&=&\frac{1}{2}|\det[z_{3}-z_{1},z_{2}-z_{1}]| ,
\text{ the area of the triangle}\nonumber\\
&=:&|\Delta_{123}|.\nonumber 
\end{eqnarray}
\end{enumerate}
\end{thm}

Again we are interested in finding two joint invariants $J_1$ and $J_2$ such that $\{J_{1,r},J_{2,r}  \}_{r=1}^k$ contains a complete fundamental set of $k$-point invariants.
We proceed similarly as for $SE(2)$ to construct a Euclidean joint invariant signature ($EJIS$ for short.) According to our general method (see Theorem \ref{star and starstar exist}), the invariants that are naturally prescribed by the result of our normalization are 
\[J_1(p_1,p_2,p_3)=|p_3-p_2| \hspace{1cm}\text{ and } \hspace{1cm} J_1(p_1,p_2,p_3)=|\Delta_{123}|.
\]
These two invariants are such that
\[  \{|p_2-p_1| \}\subset \left\{ J_1(p_i,p_{i+1},p_{i+2}), J_2(p_i,p_{i+1},p_{i+2}) \right\}_{i=k-1}^{1}
\]
with $\{|p_2-p_1| \}$ a complete set of fundamental invariants only depending  on the first two points $p_1$ and $p_2$. This guarantees the first property called two-point projectability ($\star$).

\begin{figure}[here]
\caption{}
\label{restricted domain}
\vspace{0.5cm}
\centerline{
\hbox{
\epsfysize=5.0cm
\epsfxsize=5.0cm
{\leavevmode \epsfbox{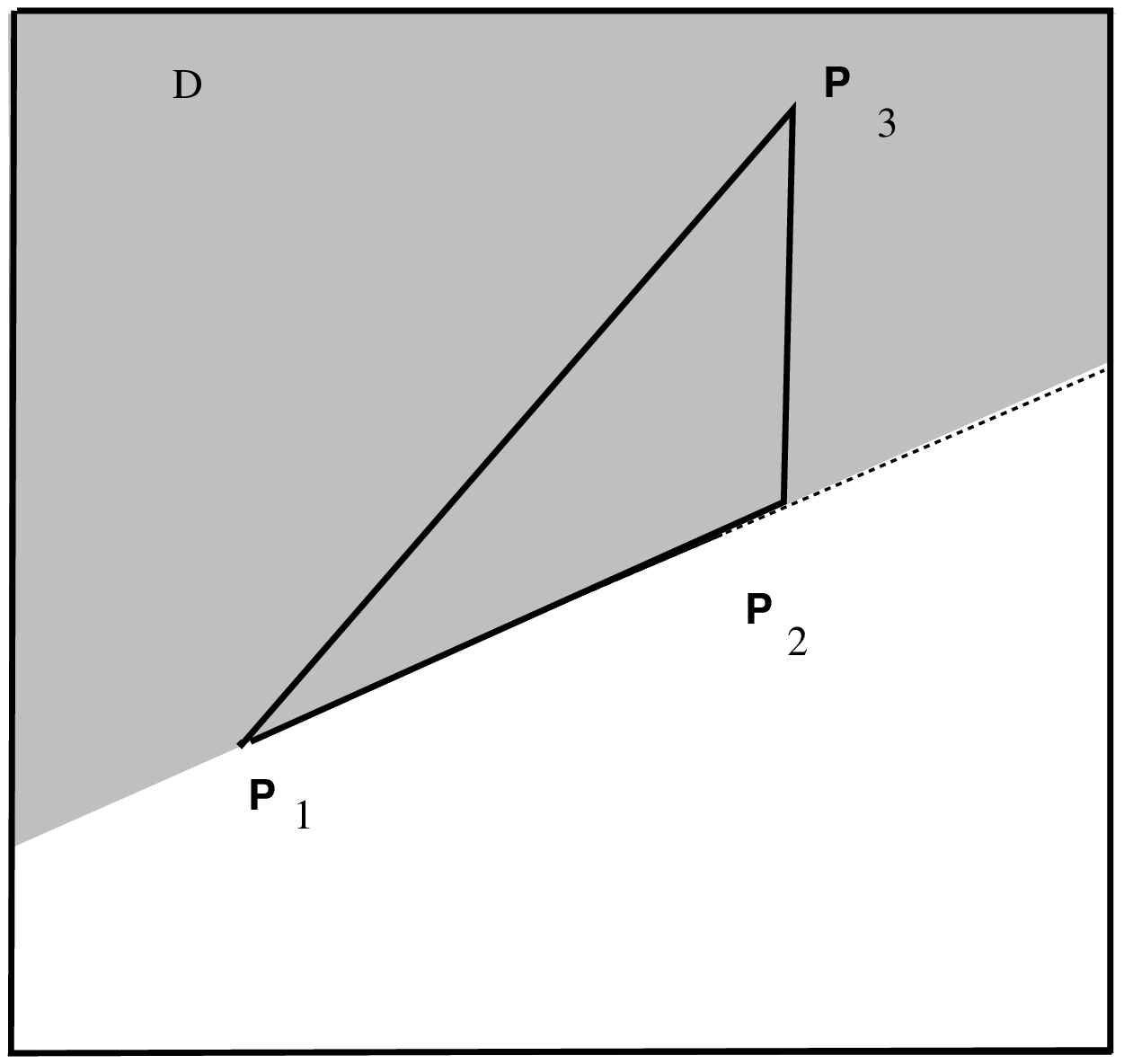}}
}
}
\end{figure}

We also have that given  $p_1, p_2, p_3\in D= \{ \Delta_{123}\geq 0 \}$, then $p_3$ is uniquely determined by the value of $J_1(p_1,p_2,p_3)$ and $J_2(p_1,p_2,p_3)$. In other words, we have 
\[ p_3=f(p_1,p_2,J_1(p_1,p_2,p_3),J_2(p_1,p_2,p_3)) 
\] 
 for $p_1,p_2,p_3\in D$. This guarantees the second property, called
third point reductivity  ($\star\star$) on the restricted domain $D$. See figure \ref{restricted domain} for an illustration.

So we can use $J_1$ and $J_2$ to parameterize a signature  for convex polygons for example, but not for all polygons. This is due to the domain restriction on  $(\star\star)$.

 Since it is desirable to be able to characterize all polygons, we would like to find a way around that difficulty. What we need is to find invariants for which $(\star\star)$ is true on a bigger domain.

Since any three-point invariant is a function of $I_1$, $I_2$ and $I_3$ and, for any $p_1$ and $p_2$, there are two choices of $p_3$ which lead to the same value of $I_1$, $I_2$ and $I_3$, there is no hope to build a signature on a bigger domain using only three-point joint invariants.  So we will try to use four-point joint invariants.

\begin{thm}\cite{Ojis}
\mylabel{second normalization E(2)}
All three-point invariants of  $E(2)$ acting on the plane are functions of
 the distances $|p_i-p_j|$, for $i,j=1,2,3$ and $i<j$.

All four-point invariants of  $E(2)$ acting on the plane are functions of
 the distances $|p_i-p_j|$, for $i,j=1,2,3,4$ and $i<j$.
\end{thm}

Observe that the fundamental three-point joint invariants written here are different than those of theorem $\ref{123 pt E(2)}$. This illustrates the non-uniqueness of  fundamental sets of invariants.

According to Theorem \ref{second normalization E(2)}, in order to have three-point projectability, it is enough that the signature contain the invariants

\[|p_2-p_1|,\quad |p_3-p_2|\quad \text{ and }\quad |p_3-p_1|.
\] 
This way, if the signature of  $(p_1,\ldots ,p_k )$ is the same as the signature of  $(q_1,\ldots, q_k )$, then we can map $(p_1,p_2,p_3)$ to $(q_1,q_2,q_3)$ with a Euclidean transformation. For example $J_1=|p_4-p_3|$ and  $J_2=|p_4-p_2|$ would do.

\begin{figure}[here]
\caption{}
\label{2 choices P4}
\vspace{0.5cm}
\centerline{
\hbox{
\epsfysize=5.0cm
\epsfxsize=5.0cm
{\leavevmode \epsfbox{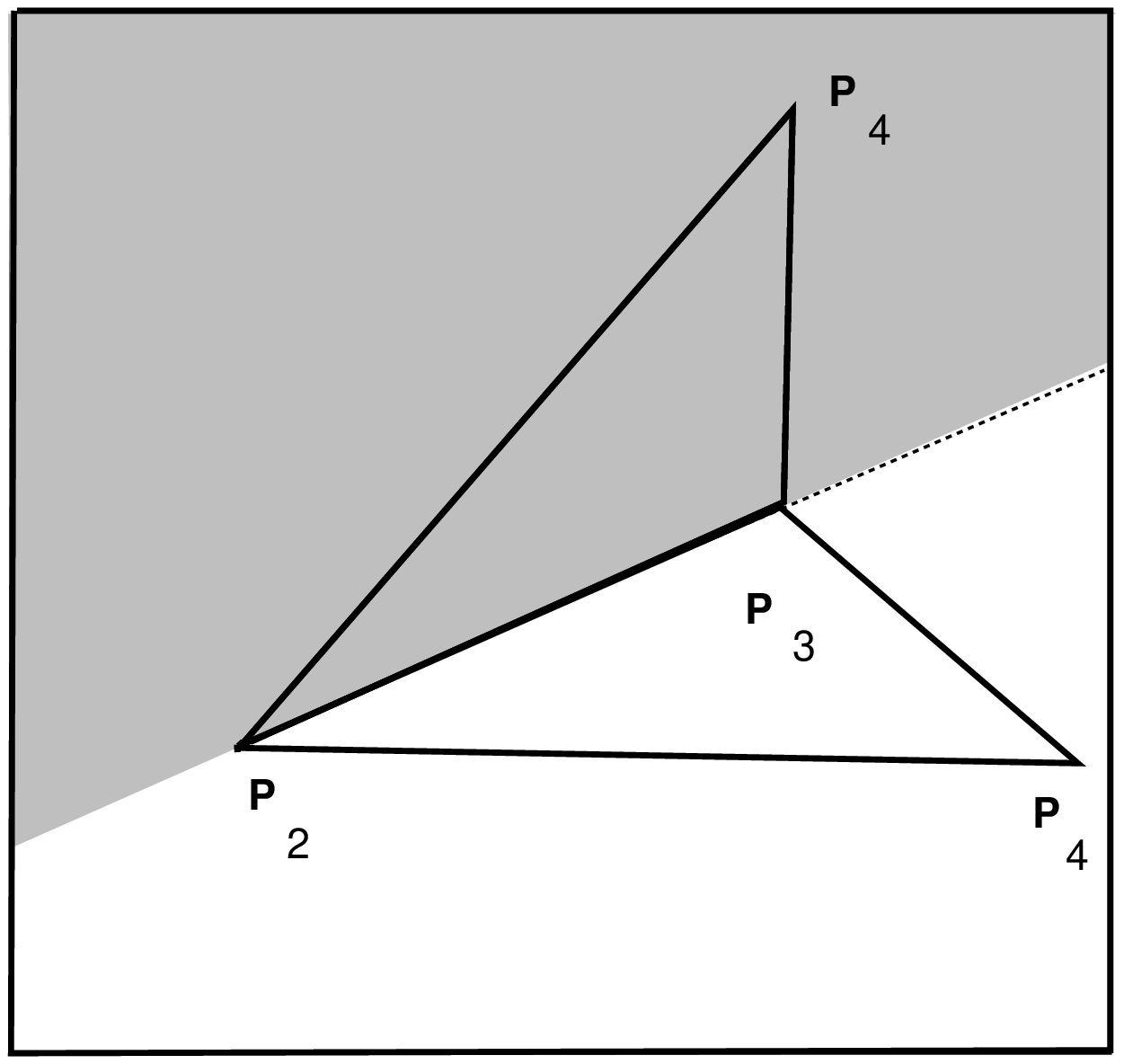}}
}
}
\end{figure}

In order to have fourth point reductivity, we need to choose two four-point joint invariants $J_1(p_1,p_2,p_3,p_4)$ and $J_2(p_1,p_2,p_3,p_4)$ which uniquely prescribe $p_4$, given $p_1$, $p_2$ and $p_3$. If we take  $J_1=|p_4-p_3|$ and $J_2=|p_4-p_2|$ then unfortunately there remains some ambiguity about the position of $p_4$ as illustrated in Figure \ref{2 choices P4}.

What we need is to know   the sign of the triangle defined by $p_2$, $p_3$ and $p_4$.
So we  look for a four-point joint  invariant  which, given  $p_1$, $p_2$ and $p_3$, determines the sign of this triangle. Observe that $\text{sign}( \Delta_{234}) $ itself is not an invariant. However  $\text{sign}(\Delta_{123}\Delta_{234})$ is an invariant (provided $\Delta_{123}\not= 0$ and $\Delta_{234}\not= 0$) and it does exactly what we need.

In fact, the invariants
\[J_1(p_1,p_2,p_3,p_4)=|p_4-p_3| \quad\text{ and }\quad J_2(p_1,p_2,p_3,p_4)=\text{sign}(\Delta_{123}\Delta_{234}) |p_4-p_2|
\] 
can be used to parameterize a Euclidean  joint invariant signature 
\[ EJIS(p_1,\ldots ,p_k)=\{J_{1,r},J_{2,r} \}_{r=1}^k\]
 for polygons for which no three consecutive points lie on a straight line. This is because
the two fundamental three-point joint invariants $|p_2-p_1|$, $|p_3-p_2|$, and  $|p_3-p_1|$ are functions
\begin{eqnarray}
|p_2-p_1|&=&f_1( \{ J_{1,r} ,J_{2,r} \}_{r=k-2}^1  ),\nonumber \\
|p_3-p_2|&=&f_2( \{ J_{1,r} ,J_{2,r} \}_{r=k-2}^1  ),\nonumber \\
|p_3-p_1|&=&f_3( \{ J_{1,r} ,J_{2,r} \}_{r=k-2}^1  ).\nonumber  
\end{eqnarray}
This guarantees property ($\star$) called {\em three-point projectability}

Moreover,  given  $p_i$, $p_{i+1}$ and $p_{i+2}$, then  $p_{i+3}$ is uniquely determined by $J_{1,i}$ and $J_{1,i}$. This guarantees property ($\star\star$) called fourth point reductivity.

We shall assume for the rest of this section that all polygons considered contain no three consecutive vertices lying on a straight line.

\begin{thm} [For polygon recognition modulo $E(2)$]
Two planar polygons $P=\langle p_1,\ldots ,p_k \rangle$ and $Q=\langle q_1,\ldots ,q_k \rangle$ are equivalent under the action of $E(2)$ 
\begin{eqnarray*}
\Leftrightarrow   EJIS(p_1,\ldots ,p_k)&=&EJIS(h\cdot (q_1,\ldots ,q_k)), \text{ for some } h\in{\mathbb H}_k,\\
\Leftrightarrow  EJIS(p_1,\ldots ,p_k)&=& c\cdot EJIS(q_1,\ldots ,q_k)\\
&  \text{ or }&\\
 EJIS(p_1,\ldots ,p_k)&=&c\cdot EJIS(q_k,\ldots ,q_1)\text{ for some } c\in{\mathbb Z}_k
\end{eqnarray*}
\end{thm}

\begin{proof}
By invariance of the functions chosen to parameterize it, the $EJIS$  of two equivalent polygons must be the same, modulo the choice of starting point and direction. This proves the necessity of the first statement. To prove necessity of the second statement, we use Lemma \ref{Hk group} and the fact that $EJIS$ commutes with rotations.

If $EJIS(p_1,\ldots ,p_k)=EJIS(h\cdot (q_1,\ldots ,q_k))$, let $(\tilde{q}_1,\ldots ,\tilde{q}_k)=h\cdot (q_1,\ldots ,q_k)$.  
Property ($\star$) allows us to conclude that $\exists g\in G$ such that  $g\cdot (p_1,p_2,p_3)=(\tilde{q}_1,\tilde{q}_2,\tilde{q}_3)$.
Property ($\star\star$) implies  that $g\cdot (p_1,\ldots, p_k)=(\tilde{q}_1,\ldots,\tilde{q}_k)$. and therefore $P\equiv Q \mod E(2)$. This proves the sufficiency of the first statement. The proof of the sufficiency of the second statement is similar.  
\end{proof}

Euclidean symmetries are of two types: rotations, which are the orientation preserving symmetries, and reflections, which are the orientation reversing symmetries. Both types of symmetries are indicated by the signature, although in general they cannot be distinguished. However for simple polygons (i.~ e.~ when its edges do not cross each other) it is possible to distinguish both types of symmetries.

\begin{thm} [For orientation preserving $E(2)$-symmetry detection in simple polygons]  
If $P=(p_1,\ldots,p_k)$ is any simple  planar polygon and $EJIS(p_1,\ldots ,p_k)=(s_1,\ldots ,s_k)$, then $P$ has a $f$-fold rotational symmetry if and only if  
\[ (s_{\frac{k}{f}+1}, \ldots ,s_k,  s_1,\ldots , s_{\frac{k}{f}})=(s_1,\ldots ,s_k),\]
that is to say, if and only if  the signature curve winds $f$ times on itself.
\end{thm}

\begin{proof}
For simple polygons, rotations are the only $E(2)$ transformations which preserve the traveling direction on the vertices, since they are the only transformations which  preserve orientation. So the proof is the same as for $SE(2)$ symmetries. 
\end{proof}

\begin{thm} [For orientation reversing $E(2)$-symmetry detection in simple polygons]

Let $P=\langle p_1,\ldots,p_k \rangle$ be any simple planar  polygon.
Then $P$ has an axe of reflection if and only if $EJIS(p_1,\ldots ,p_k)=c\cdot EJIS(p_k,\ldots ,p_1)$, for some $c\in {\mathbb Z}_k$.

More precisely,  $P$ has an axe of reflection  passing through the vertex $p_1$ if and only if
\[S(p_1,\ldots ,p_{k-1},p_k)=S(p_1,p_k,p_{k-1},\ldots ,p_2)
,\]
and $P$ has an axe of reflection  passing  in the middle of the edge joining the vertex $p_1$ to $p_2$ if and only if
\[S(p_1,p_2,\ldots ,p_{k-1},p_k)=S(p_2,p_1,p_{k},\ldots ,p_3)
.\] 
\end{thm}

\begin{proof}For simple polygons, rotations are the only $E(2)$-symmetries  which reverse the traveling direction on the vertices, since they are the only transformations which reverse orientation. By invariance of the $EJIS$ and since the $EJIS$ commutes with rotations, if $g\cdot (p_1,\ldots ,p_k)=c\cdot (p_k,\ldots ,p_1)$ for some $c\in {\mathbb Z}_k$, then $EJIS(p_1,\ldots ,p_k)=c\cdot EJIS(p_k,\ldots ,p_1)$.

Now if  $EJIS(p_1,\ldots ,p_k)=c\cdot EJIS(p_k,\ldots ,p_1)$, then ($\star$) and ($\star\star$) imply that there exists $g\in G$ such that $g\cdot (p_1,\ldots ,p_k)=c\cdot (p_k,\ldots ,p_1)$.
In particular, if $c\cdot (p_k,\ldots ,p_1)=(p_1,p_k,\ldots, p_2)$, then $p_1$ is fixed so we have an axe of reflection passing through $p_1$. Also if  $c\cdot (p_k,\ldots ,p_1)=(p_2,p_1,p_k,\ldots, p_3)$, then $p_1$ is mapped to $p_2$ and $p_2$ is mapped to $p_1$,   so we have an axe of reflection passing through the middle of the edge joining  $p_1$ to $p_2$. 
(Of course other cases can be obtained by relabeling the vertices.)  
\end{proof}



\begin{figure}[here]
\caption{}
\label{star}
\vspace{0.5cm}
\centerline{
\hbox{
\epsfysize=12.0cm
\epsfxsize=14.0cm
{\leavevmode \epsfbox{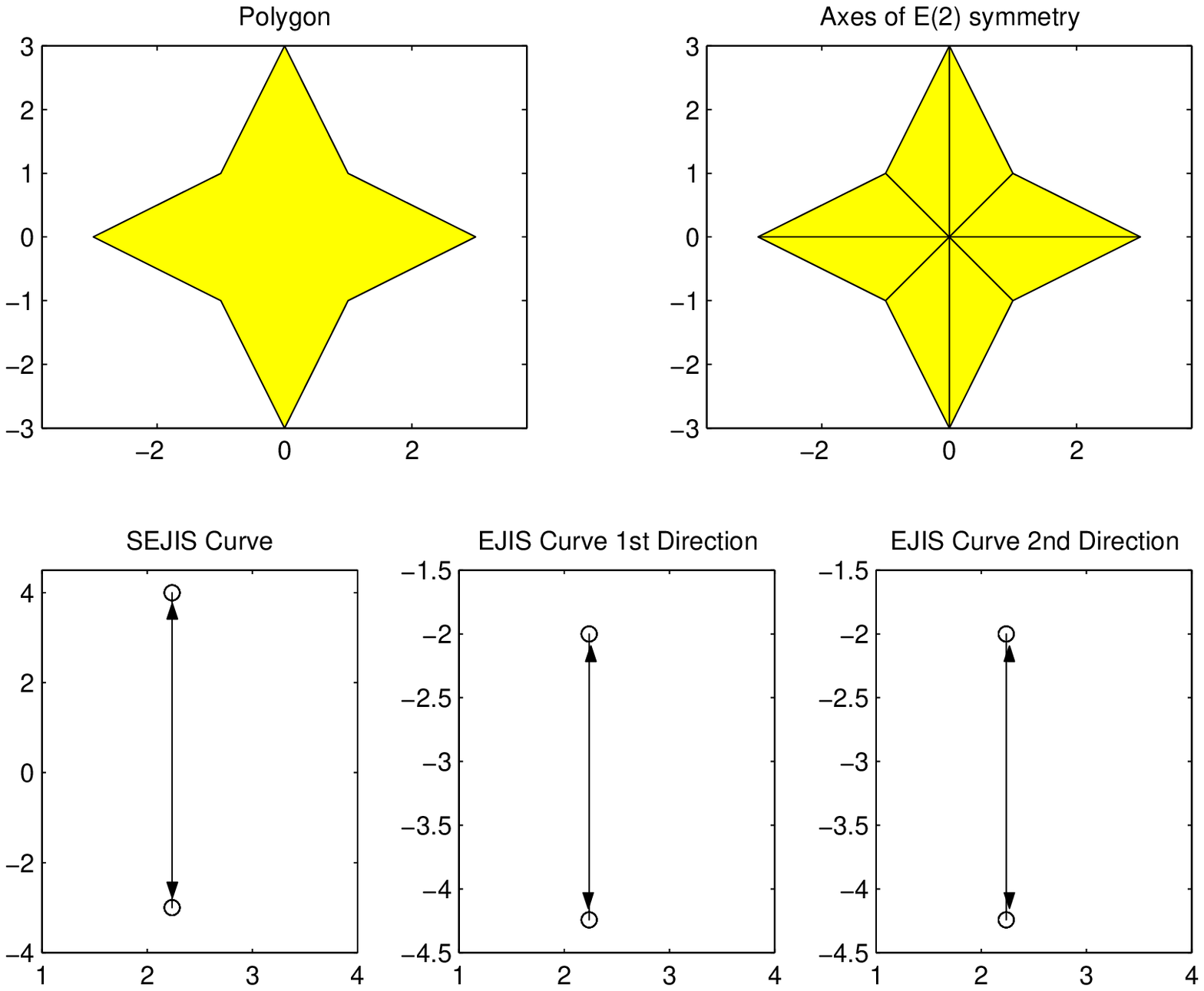}}
}
}
\end{figure}

We implemented this  algorithm using 
    Matlab and computed the results for a few examples. One of them is   the four branch star of figure \ref{star}. The program gives the following $EJIS$ (rounded to 4  digits) for the first direction we chose.

\[
EJIS1=\left[ 
\begin{array}{cccccccc} 
2.236 & 2.236 & 2.236 & 2.236 & 2.236 & 2.236 & 2.236 & 2.236  \\
-4.2426 & -2 & -4.2426 & -2 & -4.2426 & -2 & -4.2426 & -23 
\end{array}
\right]
\]

For the other direction, we obtained the following $EJIS$.

\[
EJIS2=\left[ 
\begin{array}{cccccccc} 
2.236 & 2.236 & 2.236 & 2.236 & 2.236 & 2.236 & 2.236 & 2.236  \\
 -2 & -4.2426 & -2 & -4.2426 & -2 & -4.2426 & -23  & -4.2426
\end{array} 
\right]
\]

The $EJIS$ curve, represented in figure \ref{star},  is obtained by joining those points with a straight oriented segment. Again the winding number is four, i.~ e.~ this polygon has a four fold rotational symmetry. We also detected four axes of symmetries which are also graphed on the figure.

In general we have the following.

\begin{thm}[For E(2)-symmetry detection]
If $P=\langle p_1,\ldots ,p_k \rangle $ is any planar polygon (not necessarily simple) then $P$ has a $E(2)$ symmetry if and only if there exists $h\in H_k\setminus \{ e \}$ such that   
\[EJIS(p_1,\ldots ,p_k)=EJIS(h\cdot (p_1,\ldots ,p_k))
.\] 
\end{thm}

Consider the  following instructive example.

\begin{example}Let
\[ P=\langle p_1,\ldots ,p_4 \rangle=\langle (2,1), (2,-1), (-2,1), (-2,-1)  \rangle . \] Observe that this planar polygon is not simple since two edges cross at the origin. (See Figure \ref{cross} for an illustration.)

\begin{figure}[here]
\caption{}
\label{cross}
\vspace{0.5cm}
\centerline{
\hbox{
\epsfysize=5.0cm
\epsfxsize=5.0cm
{\leavevmode \epsfbox{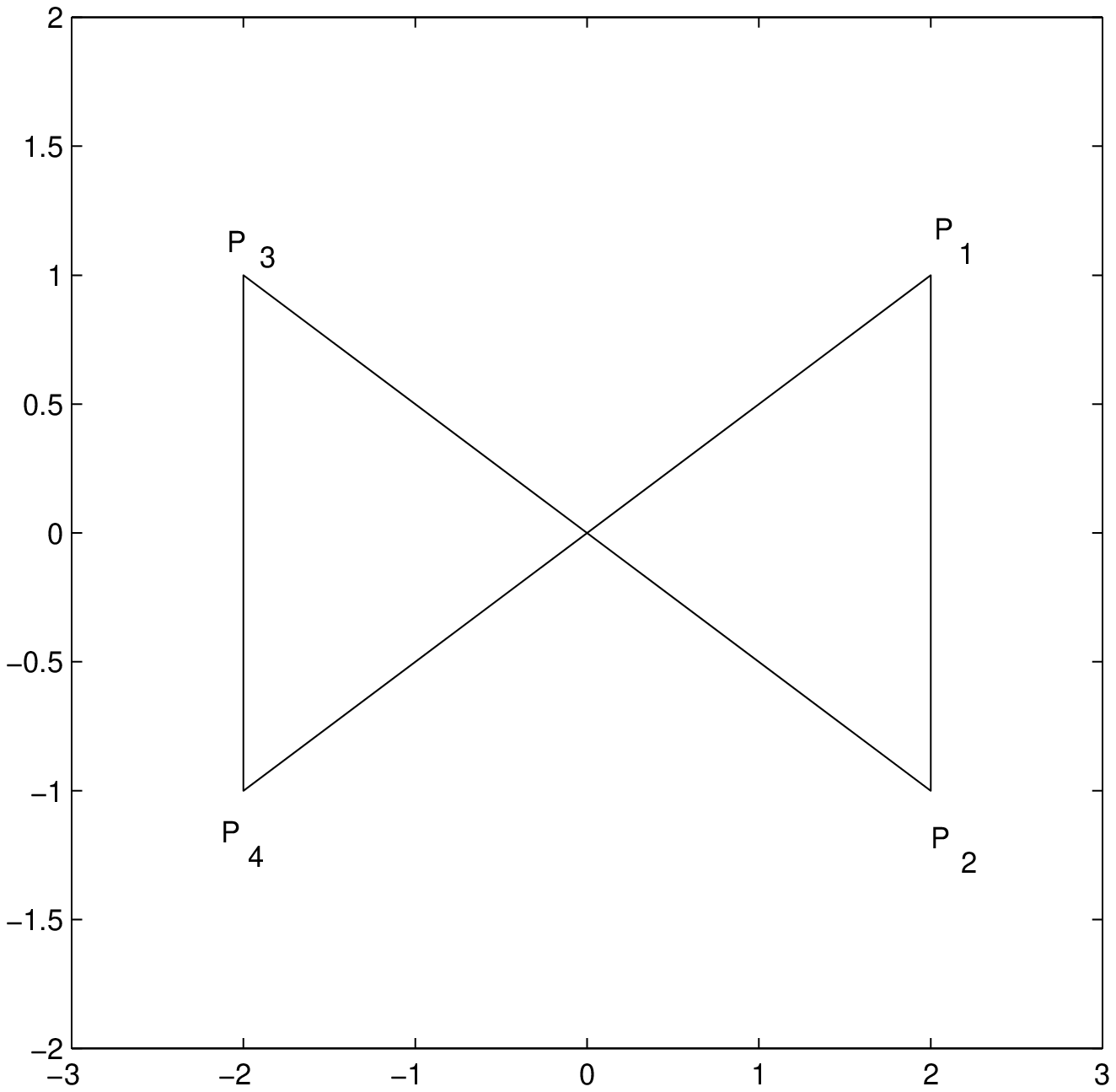}}
}
}
\end{figure}

It  has an axis of reflection which passes through the $y$-axis, and the  corresponding symmetry maps
\[(p_1,p_2,p_3, ,p_4 )\quad \text{ to }\quad (p_3,p_4,p_1,p_2  )
.\]
Although this is a orientation reversing symmetry, the traveling direction on the vertices is preserved.

Also, this polygon has a two-fold rotational symmetry which maps
\[(p_1,p_2,p_3, ,p_4 )\quad  \text{ to }\quad (p_4,p_3,p_1,p_2  )
.\]
So there is an orientation preserving symmetry which reverses the traveling direction on the vertices.
\end{example}

This example, brought to my attention by professor P.~J.~ Olver, illustrates the fact the in general, the $EJIS$ does not distinguish orientation preserving and reversing symmetries. However, using the results of the previous section, it is easy to determine which symmetries are rotations, and which symmetries are not. All one has to do is use the result provided by both the $EJIS$ and the $EJIS$  to identify which $E(2)$-symmetries are not $SE(2)$-symmetries: these are the reflections.

\section{Construction of a $G$-invariant signature curve }

Now that we have an intuitive idea of how we should build a signature curve, let us  generalize to a generic Lie group acting on a generic manifold.  The method developed will help us to construct $JIS$ curves for less intuitive Lie groups. As an illustration, the non trivial examples of the equi-affine, skewed affine and similarity groups will be presented in the last two sections.

\subsection{Two Sufficient Properties}

Our goal is to characterize $k$-gons in a $m$-dimensional manifold $M$ up to the action of an $r$-dimensional Lie group $G$.
By Theorem \ref{recognition with m invariants}, all we need to do  is to construct $m$ $n$-point joint invariants $I_1,\ldots ,I_m$ such that $\{  I_{1,r},\ldots ,I_{m,r}  \}_{r=1}^k$  contains a complete fundamental set of $k$-point joint invariants on some open set. 
One way to do this is to make sure  that they have two properties, which we call $(n-1)$-point projectability and $n^{th}$ point reductivity (see definitions below),  on some open set.  For short, we will sometimes denote them by ($\star$) and ($\star\star$) respectively. Although these conditions are stronger then needed, they have the advantage to be satisfied by the output of a simple construction algorithm.

If $n-1\leq k$, we can consider $M^{\times (n-1)}$ as a subset of $M^{\times (k)}$ by writing 
\[M^{\times (k)}=M^{\times (n-1)}\times M^{\times (k-n+1)}.
\] 
If $n-1 > k$, there exists $l\in {\mathbb N} $ such that $n-1 \leq l k$. For any such $l$, we can consider $M^{\times (n-1)}$ as a subspace of $(M^{\times (k)})^{\times (l)}$ by writing 
\[(M^{\times (k)})^{\times (l)}=M^{\times (n-1)}\times M^{\times (kl-n+1)}. 
\]

Let $l_0$ be the minimum $l\in {\mathbb N}$ such that $ n-1\leq lk$. 
Let  $c_{ir}\in {\mathbb R}$, for  $i=1,\ldots ,m$ and $ r=-n+2,\ldots ,1$. Let $C$ be the matrix $C=\{ c_{ir} \}$.
 The set
\begin{eqnarray*}
L_C:= \{ (p_1,\ldots ,p_k) \quad \text{ such that } \quad I_{i,r}(p_1,\ldots,p_n)=c_{ir},\hspace{2cm}\\
\hspace{2cm} \text{ for }\quad  i=1,\ldots ,m \quad \text{ and }\quad  r=-n+2,\ldots ,1 \}\subset M^{\times (k)}
\end{eqnarray*}
is called a {\em level set} of  $\{I_{1,r},\ldots ,I_{m,r} \}_{r=1}^k$.

We extend $L_C$  to a subset $\tilde{L}_C$ of  $(M^{\times (k)})^{\times (l_0)}$ by setting 
\[ p_{k+1}=p_1,\quad p_{k+2}=p_2,\quad \ldots\quad ,\quad p_{kl_o}=p_k .\]
The  set $\tilde{L}_C$  can  be projected in a canonical way onto a subset of $M^{\times (n-1)}$.

The first property that we will demand is the following.
\begin{definition}
We say that $m$ $n$-point joint invariants  $I_1,\ldots ,I_m: U_n\subset M^{\times (n)} \rightarrow {\mathbb R}$ are  {\em $(n-1)$-point projectable}   on $U_n$ if for any $C\in {\mathbb R}^{m\times n}$, the set $\tilde{L}_C$  
can be written as $\tilde{L}_C=U_1 \times (U_2 \cap O)$, with\\
$U_1$ an open subset of $M^{\times (kl_0-n+1)}$,\\
$U_2$ an open subset of $M^{\times (n-1)}$,\\
$O$  an orbit of $G$ acting on $M^{\times (n-1)}$.
\end{definition}

In other words,  $(n-1)$-point projectability means that  the level sets  of $\{I_{1,r},\ldots ,I_{m,r} \}_{r=1}^k$ in $(M^{\times (k)})^{\times (l_0)}$ project down to subsets of  $M^{\times (n-1)}$ which  locally correspond to orbits of the action of $G$ on $M^{(n-1)}$.  


In order to simplify the following discussion, we  introduce a new notation.
Denote by 
$ \Pi_{(i_1,i_2,\ldots ,i_R)}  U_n $ the projection 
\[  \Pi_{(i_1,i_2,\ldots ,i_R)} U_n =\{ (p_{i_1}, p_{i_2},\ldots ,p_{i_R}) | \exists (p_1,\ldots ,p_n)\in  U_n  \}\subset M^{\times (R)},
\] 
for any $1\leq i_1 < \ldots < i_R\leq n$.
Given $U_n\subset M^{\times (n)} $ and  $k\in {\mathbb N}$, let $\alpha-1=(k+n-1) \mod k$. Define  
 $\left. U_n\right|^k$ to be the set
\begin{eqnarray*}
\left. U_n\right|^k = \{(p_1,\ldots ,p_k)\in M^{\times (k)}  \text{ such that } \hspace{3cm}\\
(p_1,\ldots ,p_k,p_1,\ldots ,p_k,\quad\ldots \quad, p_1,\ldots ,p_\alpha)\in \bigcap_{i=0}^{k-1}M^{\times (i)}\times U_n\times M^{\times (k-i-1)}\}.
\end{eqnarray*}

Perhaps a more intuitive way to look at $(n-1)$-point projectability is the following.
\begin{lemma}
Let $U_n\subset M^{\times (n)}$.
The $n$-point joint invariants  $I_1,\ldots ,I_m: U_n \rightarrow {\mathbb R}$ are $(n-1)$-point projectable  on $U^{\times (n)}$ if and only if  $\{ I_{1,r},\ldots , I_{m,r} \}_{r=-n+2}^1: \left. U_n\right|^k\rightarrow {\mathbb R}$ generates a complete fundamental set of invariants on $\Pi_{(1,\ldots ,n-1)}U_n$.
\end{lemma}
\begin{proof}
By Theorem \ref{fundamental theorem}.
\end{proof}

The reason we demand $(n-1)$-point projectability is contained in this lemma.
\begin{lemma}
\mylabel{the point of star}
Let   $I_1,\ldots ,I_m: U_n \rightarrow {\mathbb R}$ be $(n-1)$-point projectable  on $U_n\subset M^{\times (n)}$ and consider the signature $S$ they define. Let $P=\langle p_1,\ldots ,p_k \rangle$ and $Q=\langle q_1,\ldots ,q_k \rangle $ be two polygons such that $(p_1,\ldots ,p_k)$ and $(q_1,\ldots ,q_k)$ are in $\left. U_n\right|^k$. If  $S(p_1,\ldots,p_k)=S(q_1,\ldots ,q_k)$, then there exists $g\in G$ such that $g\cdot (p_1,\ldots ,p_{n-1})=(q_1,\ldots ,q_{n-1})$.
\end{lemma}
\begin{proof}
By Theorem \ref{fundamental theorem}.
\end{proof}

Let $c=(c_1,\ldots ,c_m)\in {\mathbb R}^m$.
Another type of level sets are the level sets of $ I_1,\ldots ,I_m $. We will denote them by

\[\mathcal{L}_c=\{ (p_1,\ldots ,p_n)\in U_n |\quad I_i(p_1,\ldots ,p_n)=c_i, \text{ for }i=1,\ldots ,m \} \subset M^{\times (n)} .
\]

The second property we demand is the following.
\begin{definition} 
The   $n$-point joint invariants  $I_1,\ldots ,I_m: U_n\subset M^{\times (n)} \rightarrow {\mathbb R}$ are said to be  {\em $n^{th}$ point reductive  on $U_n$} if for any $c\in {\mathbb R}^m$ and any given $(p_1^0,\ldots ,p_{n-1}^0)\in \Pi_{(1,\ldots ,n-1)} U_n$,
the level set $\mathcal{L}_c$
intersects the slice
\[
\Sigma_{(p_1^0,\ldots ,p_{n-1}^0)} := \{(p_1^0,\ldots ,p_{n-1}^0, p_n)\in U_n  \}
\]
either exactly once or not at all.
\end{definition}

\begin{lemma} 
The  $n$-point joint invariants  $I_1,\ldots ,I_m: U_n\subset M^{\times (n)} \rightarrow {\mathbb R}$ 
are $n^{th}$ point reductive on $U_n$
if and only if  $p_n \in \Pi_{(n)} U_n$ is a function 
\begin{eqnarray*}
p_n&=&f(p_1,\ldots ,p_{n-1}, I_1(p_1,\ldots,p_n),\quad \ldots \quad,I_m(p_1,\ldots,p_n)\\
&=&f(p_1,\ldots ,p_{n-1}, I_{1,1},\quad \ldots \quad,I_{m,1}).
\end{eqnarray*}
\end{lemma}

\begin{definition} We say that $m$ $n$-point joint invariants are {\em perfect} on $U_n$ if they are   both $(n-1)$-point projectable ($\star$) and
 $n^{th}$ point reductive ($\star\star$) on $U_n$.
\end{definition}

\begin{prop}
\mylabel{ star and starstar implies complete}
If $I_1,\ldots ,I_m$ are perfect  on $U_n$,
then $\{ I_{1,r},\ldots , I_{m,r} \}_{r=1}^k$ contains a   complete fundamental set of  invariants on $\left. U_n\right|^k$.
\end{prop}

\begin{proof}
Let $(p_1,\ldots ,p_k)$ and $(q_1,\ldots ,q_k)\in \left. U_n\right|^k$. Consider  the signatures  $S(p_1,\ldots ,p_k)$ and  $S(q_1,\ldots ,q_k)$ parameterized by $I_1,\ldots ,I_m$. It is enough to show that   ($\star$) and   ($\star\star$) on $U_n$ imply that if $S(p_1,\ldots ,p_k)=S(q_1,\ldots ,q_k)$, then there exists $g\in G$ such that $g\cdot (p_1,\ldots , p_k)\equiv (q_1,\ldots , q_k) \mod {\mathbb H}_k$.

Assume $S( p_1,\ldots ,p_k    )= S( q_1,\ldots ,q_k )$. 
If  $n-1\geq k$, let $\beta -1= (n-1) \mod k$. Then by Lemma \ref{the point of star}, ($\star$) implies that there exists $g\in G$ such that 
\[g\cdot (p_1,\ldots ,p_k, p_1,\ldots, p_k,\quad \ldots \quad  ,p_1,\ldots ,p_{\beta})=(q_1,\ldots, q_k,q_1,\ldots ,q_k,\quad  \ldots\quad  ,q_1,\ldots ,q_{\beta}).\]
Therefore $g\cdot (p_1,\ldots ,p_k)=(q_1,\ldots ,q_k)$.

If  $n-1< k$, then  by Lemma \ref{the point of star}, ($\star$) implies that there exists $g\in G$ such that
\[g\cdot (p_1,\ldots ,p_{n-1})=(q_1,\ldots ,q_{n-1}).\]
By  ($\star\star$), we also have 
\begin{eqnarray}
g\cdot (p_2,\ldots ,p_{n})&=&(q_2,\ldots ,q_{n})\nonumber \\
g\cdot (p_3,\ldots ,p_{n+1})&=&(q_3,\ldots ,q_{n+1})\nonumber\\
&\vdots &\nonumber \\ 
g\cdot (p_{k-n+1},\ldots ,p_{k})&=&(q_{k-n+1},\ldots ,q_{k}).\nonumber
\end{eqnarray}
 Therefore $g\cdot (p_1,\ldots ,p_k)= (q_1,\ldots ,q_k)$.
\end{proof}

The converse is not true as illustrated  by the following examples. 
Take $G$ to be the special Euclidean group acting on the plane. Let $p_1$, $p_2$, and $p_3$ be three consecutive points on a polygon. Then the signed area of the triangle defined by $p_1$, $p_2$, and $p_3$ together with the distance between $p_2$, and $p_3$ satisfy ($\star$) and ($\star\star$) and therefore can be used to recognize polygons modulo orientation preserving rigid motion. However, the signature given by  the signed area of the triangle defined by $p_1$, $p_2$, and $p_3$ together with the distance between $p_1$, and $p_2$ does not satisfy   ($\star\star$), but still generates a complete fundamental set of  $k$-point invariants for any $k\in {\mathbb N}$.

Although   $n^{th}$ point reductivity is not necessary, it is an easy condition to satisfy, as will be shown later. Moreover the inverse function theorem provides an easy test for making sure this property is locally satisfied. Finally, this property is a very natural one to require when we want to detect partial equivalences in polygons.

\begin{thm}[For partial recognition modulo $G$] 
\mylabel{for partial recognition}
Let $I_1,\ldots ,I_m$ be $n$-point joint invariants which are $n^{th}$ point projectable  on $U_n\subset M^{\times (n)}$ and let $J_1,\ldots ,J_N$ be a complete fundamental set of invariants on $\Pi_{(1,\ldots ,n-1)} U_n$. Let $p_1,\ldots,p_l$ be $l$ consecutive vertices of a polygon $P$ and $q_1,\ldots,q_l$ be $l$ consecutive vertices of a polygon $Q$. Assume $l\geq n$. 
There exists $g\in G$ such that $g\cdot (p_1,\ldots,p_l)=(q_1,\ldots,q_l)$ if and only if
\begin{eqnarray}
 \label{condition 1}
J_j(p_1,\ldots ,p_{n-1})&=&J_j(q_1,\ldots ,q_{n-1}), \text{ for all } j=1,\ldots ,N \\ 
 & and &\nonumber\\ 
\label{condition 2}
I_i(p_1,\ldots ,p_n)&=&I_i(q_1,\ldots ,q_n), \\
 I_i(p_2,\ldots ,p_{n+1})&=&I_i(q_2,\ldots ,q_{n+1})\nonumber \\
&\vdots&\nonumber\\
 I_i(p_{l-n},\ldots ,p_{l})&=&I_i(q_{l-n},\ldots ,q_{l}),\text{ for all }i=1,\ldots ,m\nonumber
\end{eqnarray}

\end{thm}

\begin{proof}
By invariance of the $I$'s and $J$'s, ``$\Rightarrow$'' is true.\\
To show ``$\Leftarrow$'', assume $J_i(p_1,\ldots ,p_{n-1})=J_i(q_1,\ldots ,q_{n-1})$, for $i=1,\ldots ,N$.
 Then  there exists $g\in G$ such that  $g\cdot (p_1,\ldots ,p_{n-1})= (q_1,\ldots ,q_{n-1})$. By  ($\star\star$), condition (\ref{condition 2}) implies  that  $g\cdot (p_1,\ldots ,p_{l})= (q_1,\ldots ,q_{l})$. 
\end{proof}

If $I_1,\ldots ,I_m$ are joint invariants which are perfect, then the corresponding signature can be used for partial recognition or partial symmetry detection. 
Indeed if a complete fundamental set of invariants $J_1,\ldots ,J_N:\Pi_{(1,\ldots ,n-1)} U_n\rightarrow {\mathbb R}$ are functions 
\[J_i=f_i(\{I_{1,r},\ldots ,I_{m,r} \}_{r=-n+2}^1), \text{ for } i=1,\ldots ,N,\]
 then their value can be determined from the signature. One can therefore determine whether condition (\ref{condition 1}) is satisfied by looking at  the signatures. Condition (\ref{condition 2}) is indicated by a  partial superposition of the signatures. So both conditions can be easily checked given that we know the signatures.

\subsection{Construction of perfect $I_1,\ldots ,I_m$.}

In this section, we will determine how and in what circumstances  the moving frame method can be used to construct $m$ $n$-point joint invariants which are perfect on some open set.

Assume that the $r$-dimensional Lie group $G$ acts (locally) effectively on subsets.
Denote by $s_n$ the maximal orbit dimension of the action of $G$ on $M^{\times (n)}$.
Let  $n_o$ be the minimal integer such that for all $n\geq n_0$,   $s_n=r$, the dimension of $G$. By Theorem \ref{local freeness}, such an integer always exists.

\begin{lemma}
\mylabel{add R invariants}
Assume that  $G$ acts regularly on $U_{n+1}\subset M^{\times (n+1)}$ and $\Pi_{(1,\ldots n)}U_{n+1}$ for some $n\in {\mathbb N}$. Assume also that  $J_1,\ldots ,J_N$ is a complete fundamental set of  invariants on $\Pi_{(1,\ldots ,n)} U_{n+1}$.  Then, in a neighborhood $\tilde{U}_{n+1}$ of any point $z^{(n+1)}\in U_{n+1}$, there exist $R$  invariants $I_1,\ldots ,I_R$  such that 
\[ \{ J_1,\ldots ,J_N, I_1,\ldots ,I_R   \}
\] 
is a complete fundamental set of  invariants on $\tilde{U}_{n+1}$.
If $n\geq n_0$, then   $R=m$. Otherwise  $R<m$.
In any case, these $R$ invariants can be obtained by the moving frame normalization method or by a variant of this method. 
\end{lemma}
\begin{proof}
Theorem \ref{fundamental theorem} tells us that  there are exactly $(nm-s_n)$ fundamental $n$-point joint invariants and exactly $((n+1)m-s_{n+1})$ fundamental ($n+1$)-point joint invariants. The difference is  
\[(n+1)m-s_{n+1} - (nm-s_n)=m+s_n-s_{n+1}. \]
Let $R=m+s_n-s_{n+1}$. Observe that $R<m$, unless $s_n=s_{n+1}$.
It is shown in \cite{MBorbits} that  $s_n=s_{n+1}$ if and only if $n\geq n_0$. So if $n< n_0$, then $R$, as defined in Lemma \ref{add R invariants}, is strictly smaller than $m$, otherwise $R=m$.

Let $z^{(n+1)}\in U_{n+1}$.
Assuming that $n+1\geq n_0$, then we can build a local moving frame  $\rho(p_1,\ldots, p_{n+1})$ in a neighborhood of   $z^{(n+1)}$. Consider the group action equation   $\bar{p}_{n+1}=g\cdot p_{n+1}$. According to \cite{FO2}, setting $g=\rho(p_1,\ldots, p_{n+1})$ into this equation gives  
\[ \left. \bar{p}_{n+1}\right|_{(g=\rho(p_1,\ldots, p_n))}=(I_1,\ldots ,I_m),\]
a vector made of $m$ $(n+1)$-point invariants. Among those $m$ invariants,  there are exactly  $R$, say $I_1,\ldots I_R$, such that   $\{ J_1,\ldots ,J_N, I_1,\ldots ,I_R \}$ are functionally independent on an open subset of $U_{n+1}$. 

If $n+1 < n_0$, then a local moving frame doesn't exists in a neighborhood of  $z^{(n+1)}$. However, we can obtain a map  $\tilde{\rho}(p_1,\ldots, p_{n+1})$ such that
setting  $g=\tilde{\rho}(p_1,\ldots, p_{n+1})$ into the equation $\bar{p}_{n+1}=g\cdot p_{n+1}$ will give \[ \left. \bar{p}_{n+1}\right|_{(g=\rho(p_1,\ldots, p_n))}=(I_1,\ldots ,I_m),\]
a vector made of $m$ $(n+1)$-point invariants containing the $R$ invariants we are looking for.
\end{proof}

\begin{cor}
\mylabel{starstar does not exist}
One cannot find $m$ n-point joint invariants which are $n^{th}$ point reductive with $n \leq n_o$. 
\end{cor}

\begin{thm}
\mylabel{starstar exists}
If   $G$ acts on $U_{n+1}\subset M^{\times (n+1)}$ regularly for  some $n\geq n_0$, then in a neighborhood  of any point $z^{(n+1)}\in U_{n+1} $, there exist $m$ ($n+1$)-point joint invariants which are $(n+1)^{st}$ point reductive.
These invariants can be obtained via the moving frame method.
\end{thm}

\begin{proof}
Obtain $m$  $(n+1)$-point functionally independent invariants $\{ I_1,\ldots ,I_m\}$  as described in the proof of Lemma \ref{add R invariants}.
We claim that, on an open subset of  $U_{n+1}$, we can express $p_{n+1}$ as a function 
\[ p_{n+1}=f(p_1,\ldots ,p_n,I_1(p_1,\ldots,p_{n+1}),\ldots ,I_m(p_1,\ldots,p_{n+1})) .\] This is because if that were not the case, then the Jacobian matrix 
\[\frac{\partial (I_1,\ldots ,I_m) }{\partial (p_1,\ldots ,p_{n+1})} \]
would contain a sub-matrix
\[\frac{\partial (I_1,\ldots ,I_m) }{\partial (p_{n+1})} \] 
with rank strictly smaller than $m$, which would contradict the fact that, since $n\geq n_0$,  the invariants  $\{ I_1,\ldots ,I_m\}$ are functionally independent of invariants defined on $\Pi_{(1,\ldots ,n)} U_{n+1}$. 
\end{proof}

Let $n^\star$ be the minimum $n$ such that $G$ acts on $M^{\times (n)}$ with maximal orbit dimension $s_n < nm$. In other words, $n^\star$ is the minimum $n$ for which $n$-point joint invariants exist.

\begin{lemma}
\mylabel{translate 1}
Let $J_1,\ldots J_N$ be functionally independent invariants defined on $U_{n^\star}\subset M^{\times (n^\star)}$. 
Then 
\[\{J_{1,1},\quad \ldots \quad ,J_{N,1}, J_{1,2},\quad \ldots\quad ,J_{N,2}\}
\]
are functionally independent on an open subset of $ \left. U_{n^\star}\right|^{n^\star+1} $.
\end{lemma}
\begin{proof}
Follows from the fact that 
\[\frac{\partial (J_{1,1},\quad \ldots \quad,J_{N,1}, J_{1,2},\quad\ldots \quad,J_{N,2})}{\partial (p_1,\ldots ,p_{n^\star+1})}=\left(
\begin{array}{lr}
 M(p_1,\ldots ,p_{n^\star})&, 0_{1\times m} \\
0_{1\times m},  &M(p_2,\ldots ,p_{n^\star+1})
\end{array} 
\right)
\]
where $M$ is an $(N \times nm)$ matrix, and that,  since there are no $(n^\star-1)$-point joint invariants, the sub-matrix 
\[\frac{\partial ( J_{1,2},\quad \ldots\quad ,J_{N,2}) }{\partial (p_{n^\star+1})}
\]
has  maximal rank $m$.
\end{proof}

Let  $N_i= (n^\star+i)m-s_{n^\star+i} $ be the number of fundamental invariants of the action of $G$ on $M^{\times (n^\star+i)}$. (E.g. $N_{-1}=N_{-2}=0$.) 
We can refine the previous lemma.

\begin{lemma}
\mylabel{translate 2}
Let $n\geq n^\star$ and let $\{ J_1,\ldots J_N \}$ be a complete set of functionally independent invariants defined on 
$  U_n\subset M^{\times (n)} $.
Then there exists exactly  $(N_{n-n^\star}-N_{n-n^\star-1})$ invariants among $\{ J_{1,2},\ldots ,J_{N,2}\} $, say  $J_{1,2},\ldots , J_{(N_{n-n^\star}-N_{n-n^\star-1}),2}$, such that
\[\{J_{1,1},\quad \ldots \quad,J_{N,1}, J_{1,2},\quad \ldots\quad  , J_{(N_{n-n^\star}-N_{n-n^\star-1}),2} \}
\] 
are functionally independent  on an open subset of $ \left. U_{n^\star}\right|^{n^\star+1} $.
\end{lemma}
\begin{proof}
Follows from the fact that the rank of the sub-matrix
\[\frac{\partial (J_1,\ldots , J_N)  }{\partial p_n}
\]
is equal to $(N_{n-n^\star}-N_{n-n^\star-1})$.
\end{proof}

As a corollary, we have the following.

\begin{lemma}
\mylabel{translate 3}
Let $n\geq n^\star$ and let $J_1,\ldots J_N$ be a complete set of functionally independent invariants defined on 
$ U_n \subset M^{\times (n)} $.
Define $I_1,\ldots ,I_R$ with $R=N_{n-n^\star+1}-N_{n-n^\star}$ as in Lemma \ref{add R invariants}.
There exists
 \[ L=R-(N_{n-n^\star}-N_{n-n^\star-1})=N_{n-n^\star+1}-2N_{n-n^\star}+N_{n-n^\star-1}\]
invariants  among $I_1,\ldots ,I_R$, say $I_1\ldots ,I_{L}$, such that
\[\{ J_{1,1},\quad \ldots \quad , J_{N,1}, J_{1,2},\quad \ldots\quad , J_{N,2} ,I_1,\quad \ldots\quad ,I_{L} \}
\]
contains a complete fundamental set of invariants on an open subset of 
$ \left. U_n\right|^{n+1} $.
\end{lemma}

The proof of the next  theorem is very important as it explains the first step of the construction of perfect invariants $I_1,\ldots ,I_m$.

\begin{thm}
Let $d\in \{0,1,2,\ldots \}$ and 
let $U_{n^\star}\subset M^{\times (n^\star)}$.
If $G$ acts regularly on $\left.U_{n^\star} \right|^{n^\star+i}$,  for $i=0,1,2,\ldots ,d$, then there exist some $\bar{m}\leq m$ invariants   $I_1,\ldots ,I_{\bar{m}}$ such that
\[\{I_{1,r},\ldots ,I_{\bar{m},r} \}_{r=-n+2}^1
\]
contains a complete fundamental set of invariants on an open subset of $\left. U_{n^\star}\right|^{n^\star+d-1}$.
These can be obtained via the moving frame normalization method (or a variant of the method).
\end{thm}

\begin{proof}
By normalizing the equations  $\{ g\cdot p_i \}_{i=1}^{n^\star}$, for  $g\in G$ and $p_1,\ldots ,p_{n^\star} \in M$ as described in \cite{FO2}, we obtain functionally independent invariants  $J_1^0,\ldots ,J_{N_0}^0$ defined on an open subset of $U_{n^\star}$. 

We set
\begin{eqnarray*}
I_1 &=& J_{1,d+1}^0,\\
I_2 &=& J_{2,d+1}^0, \\
& \vdots & \\
I_{N_0} &=& J_{N_0,d+1}^0.
\end{eqnarray*}

We then normalize the equation   $g\cdot p_{n^\star+1}$ to obtain $m$  invariants $ J_1^1,\ldots ,J_{m}^1$ defined on an open subset of  $\left.U_{n^\star} \right|^{n^\star+1}$.
By Lemma \ref{add R invariants}, among those $m$ invariants, there are exactly  $R_1=N_1-N_0$, say $J_1^1,\ldots ,J_{R_1}^1$, such that
\[\{ J_1^0,\quad  \ldots\quad  ,J_{N_0}^0, J_1^1,\quad \ldots\quad ,J_{R}^1 \}\]
are functionally independent.
By Lemma \ref{translate 3}, there exists   exactly $N_1-2N_0$ invariants among $\{ J_1^1,\ldots ,J_{R_1}^1 \}$, say $J_1^1,\ldots ,J_{N_1-2N_0}^1$,  such that  
\[\{  J_{1,1}^0,\quad \ldots \quad ,J_{N_0,1}^0, J_{1,2}^0,\quad \ldots\quad ,J_{N_0,2}^0, J_1^1,\quad \ldots\quad ,J_{N_1-2N_0}^1 \}
\]
contains a complete fundamental set of invariants on an open subset of $\left.U_{n^\star} \right|^{n^\star+1}$. 

We set
\begin{eqnarray*}
 I_{N_0+1}&=&J_{1,d}^1,\\
 I_{N_0+2}&=&J_{2,d}^1,\\
& \vdots & \\
I_{N_0+N_1-2N_0}& =& J_{N_1-2N_0,d}^1.
\end{eqnarray*}
So we have now defined a  total of  $N_1-N_0$ of the $I_i$'s.

Similarly, if we normalize the equation  $g\cdot p_{n^\star+2}$ We obtain $m$ functionally independent invariants out of which  $R_2=N_2-N_1$, say $\{ J_1^2,\ldots ,J_{R_2}^2\}$, are such that
\[
J_1^0,\quad\ldots\quad ,J_{N_0}^0,
J_1^1,\quad\ldots\quad ,J_{R_1}^1 ,
J_1^2,\quad\ldots\quad ,J_{R_2}^2
\]
are functionally independent. 
By Lemma \ref{translate 3}, there exists   exactly $N_2-2N_1+N_0$ invariants among $\{ J_1^2,\ldots ,J_{R_2}^2 \}$, say $J_1^2,\ldots ,J_{N_2-2N_1+N_0}^2$,  such that
\begin{eqnarray*}
\{ J_{1,1}^0,\quad \ldots\quad  ,J_{N_0,1}^0,J_{1,2}^0,\quad \ldots \quad,J_{N_0,2}^0,J_{1,3}^0,\quad \ldots\quad  ,J_{N_0,3}^0\\
 J_{1,1}^1,\quad \ldots \quad,J_{R_1,1}^1, J_{1,2}^1,\quad \ldots\quad ,J_{R_1,2}^1, \\
J_1^2,\quad \ldots\quad ,J_{N_2-2N_1+N_0}^2 \}
\end{eqnarray*}
contains a complete fundamental set of invariants on an open subset of $\left.U_{n^\star} \right|^{n^\star+2}$ .
In other words, 
\begin{eqnarray*}
\{  J_{1,1}^0,\quad \ldots \quad,J_{N_0,1}^0, J_{1,2}^0,\quad \ldots \quad,J_{N_0,2}^0, J_{1,3}^0,\quad \ldots \quad,J_{N_0,3}^0,\\
J_{1,1}^1,\quad\ldots\quad ,J_{N_1-2N_0,1}^1, J_{1,2}^1,\quad\ldots\quad , J_{N_1-2N_0,2}^1,\\
J_1^2,\quad\ldots\quad ,J_{N_2-2N_1+N_0}^2 \}
\end{eqnarray*}
contains a complete fundamental set of invariants on an open subset of $\left. U_{n^\star}\right|^{n^\star+2}$.

We set
\begin{eqnarray*}
 I_{N_1-N_0+1}&=&J_{1,d-1}^2,\\
 I_{N_1-N_0+2}&=&J_{2,d-1}^2,\\
& \vdots & \\
I_{N_1-N_0+N_2-2N_1+N_0 }& =& J_{N_2-2N_1-N_0,d-1}^2.
\end{eqnarray*}
So we have now defined $N_2-N_1$ of the $I_i$'s.

Following this procedure $d$ times, we obtain $(N_{d-1}-N_{d-2})$ functionally independent  invariants
\[\{ I_1,\quad \ldots\quad ,I_{N_d-N_{d-1}}\},
\]
defined on some open subset of $\left. U_{n^\star}\right|^{n^\star+d}$ (although in fact they are defined on the smaller set $\Pi_{(2,\ldots ,n^\star+d)}\left. U_{n^\star}\right|^{n^\star+d}$.)
We claim that $\{ I_{1,r},\ldots ,I_{N,r} \}_{r=-n^\star+2}^1$ explicitly contains a complete fundamental set of invariants on $\Pi_{(1,\ldots ,n^\star+d-1)} \left. U_{n^\star}\right|^{n^\star+d}$. This is because, by construction, the set  
\begin{eqnarray}
\Omega =& & \{ J_{1,r}^{d-1} ,\quad \ldots \quad ,J_{N_{d-1}-2N_{d-2}+N_{d-3},r}^{d-1}  \}\nonumber \\
&\cup& \{ J_{1,r}^{d-2} ,\quad \ldots\quad  ,J_{N_{d-2}-2N_{d-3}+N_{d-4},r}^{d-2}  \}_{r=1}^2\nonumber \\
& & \vdots \nonumber \\
&\cup& \{ J_{1,r}^{0} ,\quad \ldots \quad ,J_{N_{0}-2N_{-1}+N_{-2},r}^{0}  \}_{r=1}^{d},\nonumber
\end{eqnarray}
which is contained in  $\{ I_{1,r},\ldots ,I_{N,r} \}_{r=-n^\star+2}^1$, contains exactly $N_{d-1}$ functionally independent invariants on  $\Pi_{(1,\ldots ,n^\star+d-1)} \left. U_{n^\star}\right|^{n^\star+d}$.

Observe that 
\begin{eqnarray}
N_{d-1}-N_{d-2}&=& (n^\star+d-1)m-s_{d-1} - ((n^\star+d-2)m-s_{d-2})\nonumber\\
           &=&  m-s_{d-1}+s_{d-2}\nonumber \\
           &\leq & m.\nonumber   
\end{eqnarray}
\end{proof}
Note that  if  $\{I_1,\ldots ,I_{\bar{m}} \}\subset \{I_1,\ldots ,I_m \}$, then it is guaranteed that  $\{I_1,\ldots ,I_m \}$ are  $(n^\star+d-1)$-point projectable.  
In fact, based on the   construction presented in the previous proof, we can  prove the following important result.

\begin{thm}
\mylabel{star and starstar exist}
Let $d\in \{0,1,2,\ldots   \}$ and $U_{n^\star}\subset M^{\times (n^\star)}$\\
  Assume that   $G$ acts regularly on $\left. U_{n^\star}\right|^{n^\star+i}\subset M^{\times (n^\star+i)}$, for $i=0,1,2,\ldots, d$.  
If $n^\star+d>n_0$, 
then, on an open subset of  $\left. U_{n^\star}\right|^{n^\star+d}$,  there exists $m$ ($n^\star+d$)-point joint invariants $I_1,\ldots ,I_m$  which are perfect.
These invariants can be obtained via the moving frame normalization method.
\end{thm}

\begin{proof}

We use the notation of the previous theorem.
Start with constructing $\{I_1,\ldots ,I_{N_{d-1}-N_{d-2}} \}$ as in the previous theorem. Then we normalize the equation  $g\cdot p_{n^\star+d}$, for  $g\in G$ and $p_{n^\star+d} \in M$ , and repeat the same procedure as in the previous theorem to obtain a total of $N_{d}-N_{d-1}$ invariants $I_1,\ldots ,I_{N_{d}-N_{d-1}}$. 
More precisely, we choose some $N_d-2N_{d-1}+N_{d-2}$ invariants among $J_1^d,\ldots ,J_{N_d-N_{d-1}}^d$, say $J_1^d,\ldots ,J_{N_d-2N_{d-1}+N_{d-2}}^d$, such that if $\{ J_{1,1}^{d-1},\quad \ldots \quad ,J_{N_{d-1},1}^{d-1}\}$ are a complete fundamental set of invariants defined on some open subset of $\left. U_{n^\star}\right|^{n^\star+d-1}$, then
\[J_{1,1}^{d-1},\quad\ldots\quad ,J_{N_{d-1},1}^{d-1},J_{1,2}^{d-1},\quad\ldots\quad ,J_{N_d,2}^{d-1}, J_1^d,\quad\ldots\quad ,J_{N_d-2N_{d-1}+N_{d-2}}^d
\]
are functionally independent.
We set
\begin{eqnarray*}
I_{N_{d-1}-N_{d-2}+1}&=&  J_1^d\\
I_{N_{d-1}-N_{d-2}+2}&=&  J_2^d\\
&\vdots& \\
I_{N_{d-1}-N_{d-2}+N_d-2N_{d-1}+N_{d-2}}&=&  J_{N_d-2N_{d-1}+N_{d-2}}^d, \\
\end{eqnarray*}
this obtaining $N_d-N_{d-1}$ invariants $\{ I_1,\ldots, I_{N_d-N_{d-1}}\}$.

Since  $\{I_1,\ldots ,I_{\bar{m}} \}\subset\{I_{1,r},\ldots ,I_{N_d-N_{d-1},r} \} $, the set 
\[ \{I_{1,r},\quad\ldots\quad ,I_{N_d-N_{d-1},r} \}_{r=-n^\star+2}^1\]
 explicitly contains a complete fundamental set of invariants on  $\Pi_{(1,\ldots ,n^\star+d-1)} U_{n^\star+d}$. 
Moreover, since $n^\star+d > n_0$, we have $N_d-N_{d-1}=m$.

By the same argument as in the proof of Lemma \ref{starstar exists}, we can  show that  $\{I_1,\ldots ,I_m \}$ are $(n^\star+d)^{th}$ point reductive.  
\end{proof}

An immediate corollary of Corollary \ref{starstar does not exist}, we also have

\begin{thm}
\mylabel{star and starstar do not exist}
If $n < n_0$, then for any open set $U_{n+1}\subset M^{\times (n+1)}$, there do not exist $m$ ($n+1$)-point joint invariants which 
are perfect on $U_{n+1}$.
\end{thm}

For the purpose of partial recognition, it is certainly better to use invariants depending on as few points as possible. So  one should try to build an $(n_0+1)$-point signature, which is the optimal number  for any Lie group.  However, taking more points than the minimum sometimes allows for  ($\star\star$) to be true on a bigger domain, making the detection algorithm applicable in more cases. (Recall the example of the Euclidean group acting the plane.)

The following sections contain explicit $JIS$ curves with examples  for  some slightly more difficult Lie groups namely the  the equi-affine group $SA(2)$, the skewed affine group  $SKA(2)$ and the similarity group  $SE(2)\ltimes {\mathbb R}^+$ acting on the plane.

\section{$SA(2)$ symmetry detection using  $SAJIS$  curves}

The equi-affine group $SA(2)$ is the group of area and orientation preserving  transformations in the plane. For $z\in {\mathbb R}^2$, the group transformation can be written as 
\[g\cdot z= Mz+v,\]
with $M\in SL(2)$ and  $v\in {\mathbb R}^2$. The Cartesian group action becomes  free on an open set as soon as $SA(2)$ acts on three copies of the plane. It is also regular  on $\{ (z_1,\ldots ,z_n)\in {\mathbb R}^{2})^{\times (n)}| \quad z_1,\ldots ,z_n \text{ are distinct }\}$, for all positive integers $n$. The corresponding maximal orbit dimensions are
\begin{eqnarray}
2 \text{ when } n&=&1 \nonumber \\
4 \text{ when } n&=&2 \nonumber \\
5 \text{ when } n&\geq&3 \nonumber
\end{eqnarray}

Therefore, there are no invariants of the  Cartesian action on one or two copies of the plane, while there is one fundamental invariant on three copies, and three fundamental invariants on four copies of the plane. Since $n_0+1=4$, we will try to build a four point equi-affine  joint invariant signature  ($SAJIS$.)

Let $a_{ijk}=\frac{1}{2}(z_j-z_i)\wedge(z_j-z_k)$ be the signed area of the triangle spanned by the vectors $z_j-z_i$ and $z_j-z_k$, for $z_i,z_j,z_k\in {\mathbb R}^2$.
The following  are the {\em raw} results obtained directly from the moving frame normalization method \cite{Ojis}. 

\begin{thm}

For  $SA(2)$ acting on  ${\mathbb R}^2$, we have the following. 

\begin{enumerate}

\item There are no one-point joint invariants.

\item There are no two-point joint invariants.

\item There is one fundamental three-point joint invariants
 \begin{eqnarray*}
J_1^1(p_1,p_2,p_3):({\mathbb R}^2)^{\times (3)} \rightarrow {\mathbb R},\\
 \text{namely }\quad  J_1^1(p_1,p_2)=2 a_{123}.
\end{eqnarray*}

\item There are three fundamental four-point joint invariants 
\[
J_1^2(p_1,p_2,p_3,p_4),\quad J_2^2(p_1,p_2,p_3,p_4),\quad J_3^3(p_1,p_2,p_3,p_4):({\mathbb R}^2)^{\times (4)}  \rightarrow {\mathbb R} ,
\]
namely
\begin{eqnarray}
 J_1^2(p_1,p_2,p_3,p_4)&=&2a_{123},\nonumber \\
  J_2^2(p_1,p_2,p_3,p_4)&=&-\frac{a_{134}}{a_{123}},\nonumber \\
 J_3^2(p_1,p_2,p_3,p_4)&= &2a_{124}.\nonumber
\end{eqnarray}

\end{enumerate}

\end{thm}

According to the construction described in the proof of Theorem \ref{star and starstar exist}, we take $I_1=J_{1,2}^1=2a_{234}$. Observe that $J_{1,1}^1=J_{1}^2$. So for $I_2$, we are free to take any invariant among the $J_i^2$'s {\em except} $J_{1}^2$, as long as $I_2$ and $\{ J_{1,r}^1 \}_{r=1}^2$ are functionally independent. 
In fact, we could take $I_2=-\frac{a_{134}}{a_{123}}$ or $I_2=2 a_{124}$. For simplicity, we get rid of the constants and choose to take

\[I_1 = a_{234}\quad \text{ and }\quad I_2 =  a_{124} \]
By construction, $\{I_1, I_2  \}$ are three-point projectable ($\star$) and fourth point reductive ($\star\star$) in a neighborhood of any point 
\[ z^{(3)}\in \{(p_1,p_2,p_3,p_4)\in M^{\times (3)}| \quad p_1,p_2,p_3,p_4 \text{ are distinct } \}.\]  In fact ($\star$) holds for all planar polygon assuming all four consecutive vertices are distinct.  In order to know exactly where   ($\star\star$) holds, we can solve the equations 
\[ I_1=c_1, \quad I_2=c_2 \]  
in terms for $p_4$. Computations show that a unique solution 
\[ p_4=f(p_1,p_2,p_3,I_1(p_1,p_2,p_3,p_4),I_2(p_1,p_2,p_3,p_4))\]
exists,  provided that $p_1$, $ p_2$ and $p_3$ do not lie on a straight line. Therefore our  $SAJIS$  will characterize all planar polygons for which no three consecutive vertices lie on a straight line and no consecutive vertices are identical.

We wrote a Matlab routine to test our signature on actual polygons. Using the $SAJIS$, we were able to detect  equi-affine symmetries on a collection of test polygons. One of our test polygons is shown in Figure \ref{skaffine17}. 
It is an example of a polygon with some  non-trivial affine symmetry. 
It was constructed by taking a polygon with a four-fold rotational symmetry and four axes of Euclidean symmetry and by applying a linear transformation $T\in SA(2)\setminus SE(2)$. Therefore, it has a  four-fold equi-affine symmetry and four axes of skewed-affine symmetry which are {\em not} Euclidean symmetries. 
Indeed the $SEJIS$ and $EJIS$ curves (not shown) confirmed that there is  no Euclidean symmetry. On the other hand, for a counterclockwise traveling direction, computations gave the following $SAJIS$.

\[
SAJIS=\left[ 
\begin{array}{cccccccccccc} 
1&-2&-2&1&-2&-2&1&-2&-2&1&-2&-2\\
-2&1&-2&-2&1&-2&-2&1&-2&-2&1&-2
\end{array}
\right]
\]

Our algorithm detected that, according to the $SAJIS$,  this figure has a four-fold equi-affine symmetry (winding number equal to four).

\section{$SKA(2)$ symmetry detection using  $SKAJIS$ curves}

The skewed-affine group $SKA(2)$ is the group of area  preserving transformations in the plane.   For $z\in {\mathbb R}^2$, the group action  on $z$ can be written exactly as for the previous group,  
\[g\cdot z= Mz+v,\]
where the only difference with the previous case is that  $det(M)=\pm 1$. The moving frame method gives the following.

\begin{thm}

For  $SKA(2)$ acting on  ${\mathbb R}^2$, we have the following. 

\begin{enumerate}

\item There are no one-point joint invariants.

\item There are no two-point joint invariants.

\item There is one fundamental three-point joint invariants $J_1^1(p_1,p_2,p_3):({\mathbb R}^2)^{\times (3)} \rightarrow {\mathbb R}$, namely $J_1^1(p_1,p_2)=2 |a_{123}|$.

\item There are three fundamental four-point joint invariants 
\[
J_1^2(p_1,p_2,p_3,p_4),\quad J_2^2(p_1,p_2,p_3,p_4),\quad J_3^3(p_1,p_2,p_3,p_4):({\mathbb R}^2)^{\times (4)}  \rightarrow {\mathbb R} ,
\]
namely
\begin{eqnarray}
 J_1^2(p_1,p_2,p_3,p_4)&=&2|a_{123}|,\nonumber \\
  J_2^2(p_1,p_2,p_3,p_4)&=&-\frac{a_{134}}{a_{123}},\nonumber \\
 J_3^2(p_1,p_2,p_3,p_4)&= &2|a_{124}|.\nonumber 
\end{eqnarray}
\end{enumerate}
\end{thm}

We can try to build a $n_0+1=4$ point signature with for example $I_1=|a_{234}|$ and $I_2=|a_{124}|$. But then property ($\star\star$) is only true for convex polygons. This is a very strong restriction. Inspired by our results with the Euclidean group, we choose to take 
\[I_1=sign(a_{123}a_{234})|a_{234}|\quad \text{ and }\quad I_2=sign(a_{123}a_{124}) |a_{124}|. 
\]

These two invariants satisfy property ($\star$) because $a_{123}=|I_{2,N}|$. Moreover, computations show that  we can solve for 
\[ p_4=f(p_1,p_2,p_3,I_2(p_1,\ldots ,p_4),I_2(p_1,\ldots ,p_4))\]
provided $p_1$, $p_2$ and $p_3$ do not lie on a straight line. 
So $I_1$ and $I_2$ are perfect on $\left. U_4\right|^k$, for any $k \geq 4$ and $U_4=\{(p_1,p_2,p_3,p_4)\in M^{\times (4)}|a_{123}\not=0\}$.
Therefore, $I_1$ and $I_2$ can be used to recognize all polygons $P=\langle p_1,\ldots ,p_k \rangle$ such that $(p_1,\ldots ,p_k)\in \left. U_4\right|^k$. 

\begin{figure}[here]
\caption{}
\label{skaffine17}
\vspace{0.5cm}
\centerline{
\hbox{
\epsfysize=12.0cm
\epsfxsize=14.0cm
{\leavevmode \epsfbox{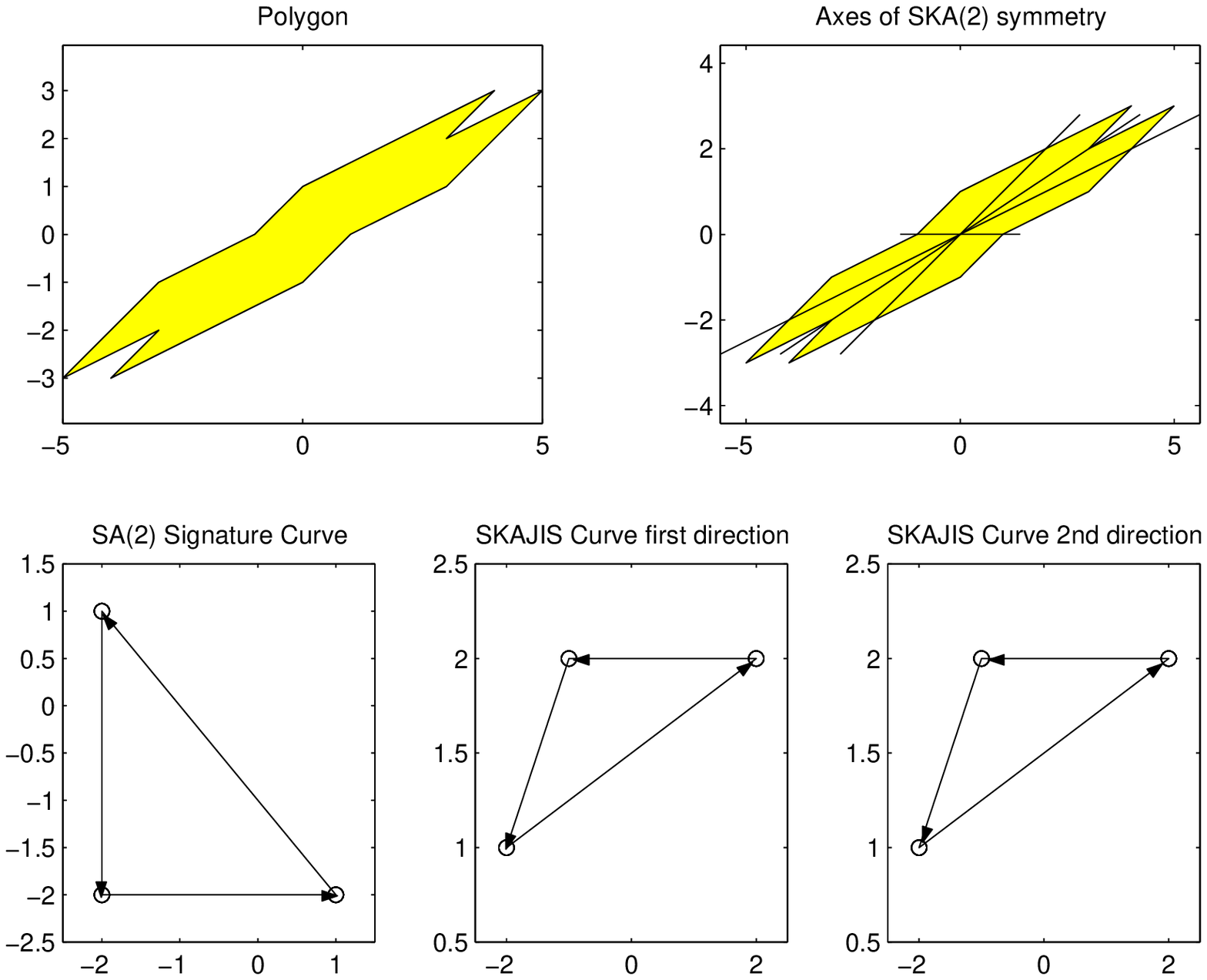}}
}
}
\end{figure}

The polygon contained in  Figure \ref{skaffine17} gives the following $SKAJIS$ for a counterclockwise orientation

\[
SKAJIS_1=\left[ 
\begin{array}{cccccccccccc} 
-2&2&-1&-2&2&-1&-2&2&-1&-2&2&-1\\
1&2&2&1&2&2&1&2&2&1&2&2
\end{array}
\right],
\]

and the following $SKAJIS$ for a clockwise direction:

\[
SKAJIS_2=\left[ 
\begin{array}{cccccccccccc} 
-2&2&-1&-2&2&-1&-2&2&-1&-2&2&-1\\
1&2&2&1&2&2&1&2&2&1&2&2
\end{array}
\right].
\]

Our algorithm detected that, according to the $SKAJIS$,  this figure has a four-fold equi-affine symmetry (winding number equal to four)  and four axes of skewed-affine symmetry (shown in figure \ref{skaffine17}.) 

\section{$SIM(2)$ symmetry detection using  $SIMJIS$ curves}

\begin{figure}[here]
\caption{}
\label{simsejis2223}
\vspace{0.5cm}
\centerline{
\hbox{
\epsfysize=12.0cm
\epsfxsize=14.0cm
{\leavevmode \epsfbox{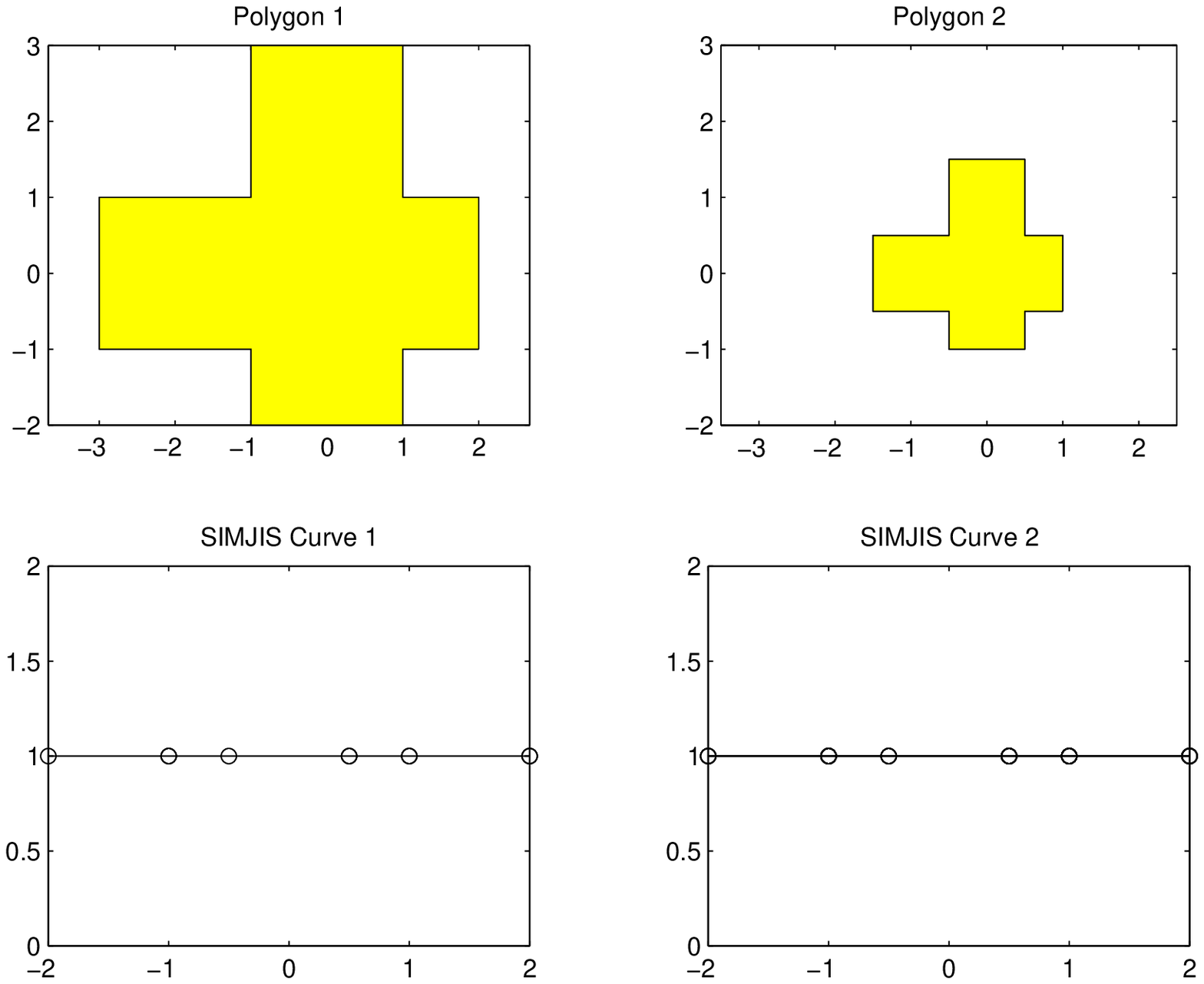}}
}
}
\end{figure}

Another important group is the similarity group given by all  special Euclidean and scaling transformations of the form
\[\bar{x}=\lambda M\cdot x+b
\]
with $\lambda\in {\mathbb R}^+$, $M\in SL(2)$ and $b\in {\mathbb R}^2$. Observe that this group acts locally freely and transitively on $\{(x_1,x_2)\in ({\mathbb R}^2)^{\times (2)}|\quad x_1\not= x_2 \}$. Using the moving frame method and following the construction described in this paper, we chose to take the following two invariants
\[I_1(p_1,p_2,p_3)=\frac{a_{123}}{|p_2-p_1|^2}, \hspace{1cm}I_1(p_1,p_2,p_3)=\frac{(p_2-p_1)\cdot (p_3-p_1)}{|p_2-p_1|^2},
\]
which are perfect on $\{(x_1,x_2,x_3)\in  ({\mathbb R}^2)^{\times (3)}|\quad x_1,x_2,x_3 \text{ are distinct }\}$. 
In order to test the corresponding similarity joint invariant signature ($SIMJIS$), we computed the $SIMJIS$ associated to a collection of polygons with some rotational symmetry and checked that the signature did illustrate the symmetry.  Observe that a polygon {\em cannot} have a scaling symmetry. So only rotational symmetries are indicated by the $SIMJIS$. The scaling part of the similarity group is of interest when comparing two polygons.  An example of two polygons equivalents under a scaling transformation is presented in Figure \ref{simsejis2223}.  The associated signature for a counter clockwise orientation
\[SIMJIS=\left[ 
\begin{array}{cccccccccccc} 
0.5&-2&1&1&-1&1&1&-0.5&2&0.5&-1&2\\
1&1&1&1&1&1&1&1&1&1&1&1
\end{array}
\right],
\]
is the same in both cases. For clarity, we did not graph the arrows representing the direction of each segment joining  consecutive points of the signature curve.

\section*{Acknowledgments}
I want to thank my advisor Peter J.~ Olver for his advice and support. I also want to thank Irina Kogan for her friendly encouragements and many  useful comments.

\newpage
\bibliography{bib}
\bibliographystyle{abbrv}
\end{document}